\documentclass{amsproc}

\numberwithin{equation}{section}

\theoremstyle{plain}
\newtheorem{introthm}{Theorem}   

\newtheorem{thm}[equation]{Theorem}  
\newtheorem{cor}[equation]{Corollary}     
\newtheorem{lem}[equation]{Lemma}         
\newtheorem{prop}[equation]{Proposition} 
\newtheorem{addendum}[equation]{Addendum}

\theoremstyle{definition}

\theoremstyle{remark}

\usepackage{amssymb}
\usepackage{amsmath}
\usepackage[matrix,arrow,curve,cmtip]{xy}

\begin{document}

\title{Bi-relative algebraic $K$-theory and\\ topological cyclic homology} 

\author{Thomas Geisser}

\address{University of Southern California, Los Angeles, California}

\email{geisser@math.usc.edu}

\author{Lars Hesselholt}

\address{Massachusetts Institute of Technology, Cambridge,
Massachusetts}

\email{larsh@math.mit.edu}

\thanks{Both authors were supported in part by the National Science
  Foundation (USA). The second named author received additional
  support from the COE (Japan).}

\maketitle

\section*{Introduction}

It is well-known that algebraic $K$-theory preserves products of
rings. However, in general, algebraic $K$-theory does not preserve
fiber products of rings, and one defines bi-relative algebraic
$K$-theory so as to measure the deviation. It was proved recently by
Corti\~{n}as~\cite{cortinas} that, rationally, bi-relative algebraic
$K$-theory and bi-relative negative cyclic homology agree. In this
paper, we show that, with finite coefficients, bi-relative algebraic
$K$-theory and bi-relative topological cyclic homology agree. As an
application, we show that for a possibly singular curve over a perfect
field $k$ of positive characteristic $p$, the cyclotomic trace map
induces an isomorphism of the $p$-adic algebraic $K$-groups and the
$p$-adic topological cyclic homology groups in non-negative
degrees. As a further application, we show that the difference between
the $p$-adic $K$-groups of the integral group ring of a finite group
and the $p$-adic $K$-groups of a maximal $\mathbb{Z}$-order in the
rational group algebra can be expressed entirely in terms of
topological cyclic homology.

Let $F$ be a functor that to a unital associative ring $A$ associates
a symmetric spectrum $F(A)$. If $I \subset A$ is a two-sided ideal,
the relative term $F(A,I)$ is defined to be the homotopy fiber of the
map $F(A) \to F(A/I)$ induced from the canonical projection. If
further $f \colon A \to B$ is a ring homomorphism such that $f \colon
I \to f(I)$ is an isomorphism onto a two-sided ideal of $B$, there is
an induced map $F(A,I) \to F(B,f(I))$, and the bi-relative term
$F(A,B,I)$ is defined to be the homotopy fiber of this map. Then there
is a distinguished triangle in the stable homotopy category
$$F(A,B,I) \to F(A,I) \xrightarrow{f_*} F(B,f(I))
\xrightarrow{\partial} \Sigma F(A,B,I)$$
which is natural with respect to the triple $(A,B,I)$. We write the
induced long-exact sequence of homotopy groups as
$$\cdots \to F_q(A,B,I) \to F_q(A,I) \xrightarrow{f_*} F_q(B,f(I))
\xrightarrow{\partial} F_{q-1}(A,B,I) \to \cdots.$$
The theorem of Corti\~{n}as~\cite{cortinas} states that the trace map
$$K_q(A) \otimes \mathbb{Q} \to
\operatorname{HC}_q^{-}(A \otimes \mathbb{Q})$$ 
from rational $K$-theory of $A$ to negative cyclic homology of
$A \otimes \mathbb{Q}$ induces an isomorphism of the associated
bi-relative theories. In this paper we prove the following analogous
result for the $K$-theory with finite coefficients:

\begin{introthm}\label{main}Let $f \colon A \to B$ be a map of unital
associative rings, let $I \subset A$ be a two-sided ideal and assume
that $f \colon I \to f(I)$ is an isomorphism onto a two-sided ideal of
$B$. Then the map induced by the cyclotomic trace map
$$K_q(A,B,I,\mathbb{Z}/p^v) \to
\{\operatorname{TC}_q^n(A,B,I;p,\mathbb{Z}/p^v)\}_{n \geqslant 1}$$
is an isomorphism of pro-abelian groups, for all integers $q$, all
primes $p$, and all positive integers $v$. 
\end{introthm}

We recall from Bass~\cite[Thm.~XII.8.3]{bass} that the common
pro-abelian group is zero, for $q \leqslant 0$. We also remark that
Thm.~\ref{main} implies the slightly weaker statement that the map of
spectra induced by the cyclotomic trace map
$$K(A,B,I) \to \operatorname{TC}(A,B,I;p)
= \operatornamewithlimits{holim}_n \operatorname{TC}^n(A,B,I;p)$$
becomes a weak equivalence after $p$-completion. However, for the
proof, it is essential that we work with the pro-abelian groups and
not pass to the limit. We also mention that Dundas and
Kittang~\cite{dundaskittang} have shown that Thm.~\ref{main} implies
the analogous statement for connective symmetric ring spectra.

A weaker version of the theorem of Corti\~{n}as was first stated as a
conjecture by Geller, Reid, and
Weibel~\cite{gellerreidweibel,gellerreidweibel1} and used to evaluate
the rational algebraic $K$-theory of curves over a field of
characteristic zero. In a similar manner, we obtain the following
result from Thm.~\ref{main}:

\begin{introthm}\label{curves}Let $X \to \operatorname{Spec} k$ be any
curve over a field $k$ of positive characteristic $p$. Then, for all
positive integers $v$, the cyclotomic trace map
$$K_q(X,\mathbb{Z}/p^v) \to
\{\operatorname{TC}_q^n(X;p,\mathbb{Z}/p^v)\}_{n \geqslant 1}$$
is an isomorphism of pro-abelian groups in degrees $q \geqslant r$,
where $p^r = [k : k^p]$ is the degree of $k$ over the subfield
$k^p$ of $p$th powers.
\end{introthm}

In the proof of Thm.~\ref{curves}, which is given in
Sect.~\ref{proofoftheorem} below, we use Thm.~\ref{main} to reduce to
the smooth case which we have proved earlier~\cite[Thm.~4.2.2]{gh}. 
We remark again that Thm.~\ref{curves} implies the slightly weaker
statement that the map of spectra induced by the cyclotomic trace map
$$K(X) \to \operatorname{TC}(X;p)
= \operatornamewithlimits{holim}_{n} \operatorname{TC}^n(X;p)$$
becomes a weak equivalence after $p$-completion and after passage
to $(r-1)$-connected covers. The topological cyclic homology groups on
the right-hand side of the statement of Thm.~\ref{curves} can often be
effectively calculated. As an example, we have the following
result; compare~\cite{gellerreidweibel}. The proof is similar to the
proof of~\cite[Thm.~E]{hm} and~\cite[Thm.~A]{hm2} and will appear
in~\cite{h4}.

\begin{introthm}\label{coordinateaxes}Let $k$ be a regular
$\mathbb{F}_p$-algebra, let $A = k[x,y]/(xy)$ be the coordinate ring
of the coordinate axes in the affine $k$-plane, and let $I \subset A$
be the ideal generated by the variables $x$ and $y$. Then for all
positive integers $q$, there is a canonical isomorphism
$$K_q(A,I) \xleftarrow{\sim}
\bigoplus_{m \geqslant 1} \mathbf{W}_m\Omega_k^{q-2m}$$
where $\mathbf{W}_m\Omega_k^j$ is the group of big de~Rham-Witt
$j$-forms of $k$~\cite{hm2,h3}.
\end{introthm}

The group $K_2(A,I)$ was evaluated twenty-five years ago by Dennis and
Krusemeyer~\cite{denniskrusemeyer}, but it was previously only known
that the higher relative $K$-groups are $p$-primary torsion
groups~\cite{weibel2}. We point out that, at present, the structure of
the $K$-groups of singular schemes is very far from understood. In
particular, there is currently no motivic cohomology theory that
explain the structure of the $K$-groups that we find here. To 
illustrate this, let $k$ be a perfect field of positive characteristic
$p$. Then Thm.~\ref{coordinateaxes} shows that the $p$-adic $K$-groups
of $X = \operatorname{Spec} k[x,y]/(xy)$ are concentrated in even
degrees with the group in positive degree $2m$ canonically isomorphic
to the group $\mathbf{W}_m(k)$ of big Witt vectors of length $m$ in
$k$. In particular, the $p$-adic $K$-groups of $X$ cannot be generated
from $K_1(X)$ by means of products and transfers as is the case for
smooth schemes.

Finally, let $A = \mathbb{Z}[G]$ be the integral group ring of a
finite group, and let $B = \mathfrak{M}$ be a maximal
$\mathbb{Z}$-order of the rational group algebra. Then Thm.~\ref{main}
shows that, for every prime $p$, the difference between the $p$-adic
$K$-groups of $A$ and $B$ can be expressed entirely in terms of
topological cyclic homology. Indeed, if we choose $B$ such that $A
\subset B$, then the ideal $I = |G| \cdot B$ is a common ideal of both
$A$ and $B$, and the quotient rings $A/I$ and $B/I$ are finite. The
precise relationship is most clearly expressed for $G$ a finite
$p$-group where we obtain the following result. We are grateful to
Holger Reich for drawing our attention to this application of
Thm.~\ref{main}.

\begin{introthm}\label{grouprings}Let $A = \mathbb{Z}[G]$ be the
integral group ring of a finite $p$-group, and let $B = \mathfrak{M}$
be a maximal $\mathbb{Z}$-order of the rational group algebra such
that $A \subset B$. Then there is a natural sequence of pro-abelian
groups
$$\cdots \to
K_q(A,\mathbb{Z}/p^v) \to
\begin{Bmatrix}
{K_q(B,\mathbb{Z}/p^v)} \cr
{\oplus} \cr
{\operatorname{TC}_q^n(A;p,\mathbb{Z}/p^v)} \cr
\end{Bmatrix} \to
\{\operatorname{TC}_q^n(B;p,\mathbb{Z}/p^v)\} \to \cdots$$
which is exact, for all non-negative integers $q$, and all positive
integers $v$.
\end{introthm}

We again remark that Thm.~\ref{grouprings} implies the slightly weaker
result that the following square of spectra becomes homotopy-cartesian
after $p$-completion and passage to connective covers:
$$\xymatrix{
{K(\mathbb{Z}[G])} \ar[r] \ar[d] &
{\operatorname{TC}(\mathbb{Z}[G];p)} \ar[d] \cr
{K(\mathfrak{M})} \ar[r] &
{\operatorname{TC}(\mathfrak{M};p).} \cr }$$
The strength of Thm.~\ref{grouprings} lies in the fact that the
maximal order $\mathfrak{M}$ is a regular ring whose $K$-groups
fit in the following localization sequence~\cite[Thm.~1.17]{oliver}:
$$\cdots \to \bigoplus_{\ell} K_q(\mathfrak{M}_{\ell}/J_{\ell}) \to
K_q(\mathfrak{M}) \to
K_q(\mathbb{Q}[G]) \to \cdots.$$
The sum on the left-hand side ranges over all rational primes $\ell$
and $J_{\ell} \subset \mathfrak{M}_{\ell}$ is the Jacobson radical in
the $\ell$-adic completion of $\mathfrak{M}$. The rings
$\mathbb{Q}[G]$ and $\mathfrak{M}_{\ell}/J_{\ell}$ are semi-simple.

The proof of Thm.~\ref{main} is similar, in outline, to the proof of
the theorem of Corti\~{n}as~\cite{cortinas} and uses two key
observations by Cuntz and Quillen~\cite{cuntzquillen}. The first
observation is that it suffices to prove Thm.~\ref{main} in the case
where the non-unital associative ring $I$ can be embedded as a
two-sided ideal of a free unital associative ring. The second
observation is that, in this case, the pro-algebra $\{I^m\}$ is
homological unital in the sense that for all primes $p$, and all
positive integers $q$, the following pro-abelian group is zero:
$$\{ \operatorname{Tor}_q^{\mathbb{Z} \ltimes
  I^m}(\mathbb{Z},\mathbb{Z}/p) \}_{m \geqslant 1}.$$
Here $\mathbb{Z} \ltimes I^m$ is the unital associative ring obtained
from $I^m$ by adjoining a unit. In Sect.~\ref{ktheorysection} below
we show that this implies that for all integers $q$, all primes $p$,
and all positive integers $v$, the pro-abelian group
$$\{K_q(A,B,I^m,\mathbb{Z}/p^v)\}_{m \geqslant 1}$$
is zero. This result is a generalization of the excision theorem of
Suslin-Wodzicki~\cite{suslin4,suslinwodzicki}. Similarly, we show in
Sect.~\ref{thhsection} that for all integers $q$, all primes $p$, and
all positive integers $v$, the pro-abelian group
$$\{\operatorname{TC}_q^n(A,B,I^m;p,\mathbb{Z}/p^v)\}_{m,n \geqslant 1}$$ 
is zero. Hence, if $I$ can be embedded as an ideal of a free unital
associative ring, the cyclotomic trace map induces an isomorphism of
pro-abelian groups in the top row of the following diagram:
$$\xymatrix{
{\{K_q(A,B,I^m,\mathbb{Z}/p^v)\}_{m \geqslant 1}}
\ar[r]^(.455){\sim} \ar[d] & 
{\{\operatorname{TC}_q^n(A,B,I^m;p,\mathbb{Z}/p^v)\}_{m,n \geqslant 1}}
\ar[d] \cr
{K_q(A,B,I,\mathbb{Z}/p^v)} \ar[r] &
{\{\operatorname{TC}_q^n(A,B,I;p,\mathbb{Z}/p^v)\}_{n \geqslant 1}.} \cr
}$$
Finally, it follows from a theorem of McCarthy~\cite{mccarthy1} and
from~\cite[Thm.~2.2.1]{gh3} that also the lower horizontal map is an
isomorphism of pro-abelian groups. The details of this argument are
given in Sect.~\ref{proofoftheorem} below. This completes the
outline of the proof of Thm.~\ref{main}.

A pro-object of a category $\mathcal{C}$ is a functor from a directed
partially ordered set to the category $\mathcal{C}$, and a
\emph{strict} map between two pro-objects with the same indexing set
is a natural transformation. A general map from a pro-object $X =
\{X_i\}_{i \in I}$ to a pro-object $Y = \{Y_j\}_{j \in J}$ is an
element of the set
$$\operatorname{Hom}_{\operatorname{pro}-\mathcal{C}}(X,Y)=
\operatornamewithlimits{lim}_J
\operatornamewithlimits{colim}_I
\operatorname{Hom}_{\mathcal{C}}(X_i,Y_j).$$
In particular, a pro-object $X = \{X_i\}_{i \in I}$ in a category with
a null-object is zero if for all $i \in I$, there exists
$i' \geqslant i$ such that the map $X_{i'} \to X_i$ is zero. We write
$X[-]$ for a simplicial object in a category $\mathcal{C}$, and we
write $X[k]$ for the object in simplicial degree $k$. If $I$ is a
non-unital associative ring, we write $\mathbb{Z} \ltimes I$ for the
unital associative ring given by the product $\mathbb{Z} \times I$
with multiplication $(n,x) \cdot (n',x') = (nn',nx' + n'x + xx')$. 

This paper was written in part while the authors visited the
University of Tokyo. We would like to express our sincere gratitude to
the university and to Takeshi Saito for their kind hospitality and for
the stimulating atmosphere they provided. Finally, it is a particular
pleasure to acknowledge the help we have received from Guillermo
Corti\~{n}as. He long ago suggested to us to prove the main theorem of
this paper by adapting the integral excision theorem of
Suslin~\cite{suslin4} to the pro-setting. He also pointed out a gap in
our arguments in an earlier version of this paper.

\section{$K$-theory}\label{ktheorysection}

In this section we prove the following variant of the
excision theorem of Suslin-Wodzicki~\cite{suslinwodzicki,suslin4}.
Our treatment closely follows Suslin~\cite{suslin4}.

\begin{thm}\label{ktheory}Let $f \colon A \to B$ be a map of
unital associative rings, and let $I \subset A$ be a two-sided ideal
such that $f \colon I \to f(I)$ is an isomorphism onto a two-sided
ideal of $B$. Let $p$ be a prime, and suppose that for all positive
integers $q$, the pro-abelian group
$$\{ \operatorname{Tor}_q^{\mathbb{Z} \ltimes I^m}
(\mathbb{Z},\mathbb{Z}/p) \}_{m \geqslant 1}$$
is zero. Then for all integers $q$, and all positive integers $v$, the
map of pro-abelian groups
$$\{ K_q(A,I^m,\mathbb{Z}/p^v) \}_{m \geqslant 1} \to
\{ K_q(B,f(I)^m,\mathbb{Z}/p^v) \}_{m \geqslant 1}$$
is an isomorphism.
\end{thm}

We will use a \emph{functorial} version of the plus-construction of
Quillen. So
$$\lambda \colon X \to X^+$$
will be a natural transformation of functors from the category of
pointed spaces and base-point preserving continuous maps to itself.
For example, one can use the functorial
localization~\cite[Def.~1.2.2]{hirschhorn} with respect to a pointed
space $W$ whose reduced integral homology is trivial and whose
fundamental group admits a non-trivial group homomorphism to
$GL(\mathbb{Z})$~\cite{berrickdwyer}.

Let $A$ be a unital associative ring, and let $I \subset A$ be a
two-sided ideal. Let $GL(I)$ and $\overline{GL}(A/I)$ be the kernel
and image of the map $GL(A) \to GL(A/I)$ induced by the canonical
projection. We let $F(A,I)$ be the homotopy fiber of the canonical
map $BGL(A)^+ \to B\overline{GL}(A/I)^+$, and let $X(A,I)$ be the
homotopy pull-back of $B\overline{GL}(A/I)$ and $BGL(A)^+$ over
$B\overline{GL}(A/I)^+$. We then have the following diagram of pointed
spaces and base-point preserving continuous maps where the rows are
fiber sequences:
$$\xymatrix{
{BGL(I)} \ar[r] \ar[d] &
{BGL(A)} \ar[r] \ar[d] &
{B\overline{GL}(A/I)} \ar@{=}[d] \cr
{F(A,I)} \ar[r] \ar@{=}[d] &
{X(A,I)} \ar[r] \ar[d] &
{B\overline{GL}(A/I)} \ar[d] \cr
{F(A,I)} \ar[r] &
{BGL(A)^+} \ar[r] &
{B\overline{GL}(A/I)^+.} \cr
}$$
We consider the relative Serre spectral sequence
$$\begin{aligned}
E_{s,t}^2 & = H_s(\overline{GL}(A/I),H_t(F(A,I),BGL(I),\mathbb{Z}/p)) \cr
{} & \Rightarrow H_{s+t}(X(A,I),BGL(A),\mathbb{Z}/p) \cr
\end{aligned}$$
which converges to zero. Indeed, both of the maps $BGL(A) \to
BGL(A)^+$ and $X(A,I) \to BGL(A)^+$ induce isomorphisms on homology,
and hence so does the map $BGL(A) \to X(A,I)$. Hence, the relative
homology groups, defined as the reduced homology groups of the mapping
cone, are zero. The diagram and the spectral sequence considered above
both are functorial in the pair $I \subset A$. Hence, if we consider
the pro-ideal $\{I^m\}_{m \geqslant 1}$ of powers of the ideal $I
\subset A$, we obtain a spectral sequence of pro-abelian groups.

\begin{lem}\label{inductionlemma}If the pro-abelian group
$\{H_t(F(A,I^m),BGL(I^m),\mathbb{Z}/p)\}$ is zero, for
$0 \leqslant t < q$, then the pro-abelian group 
$$\{H_s(\overline{GL}(A/I^m),H_q(F(A,I^m),BGL(I^m)),\mathbb{Z}/p)\}$$
is zero, for $0 \leqslant s \leqslant 1$. If, in addition, the
canonical map
$$\{H_{q-1}(BGL(I^m),\mathbb{Z}/p)\} \to
\{H_0(\overline{GL}(A/I^m),H_{q-1}(BGL(I^m),\mathbb{Z}/p))\}$$
is an isomorphism of pro-abelian groups, then the pro-abelian group
$$\{H_q(F(A,I^m),BGL(I^m),\mathbb{Z}/p)\}$$
is zero.
\end{lem}

\begin{proof}We omit the $\mathbb{Z}/p$-coefficients for the homology
throughout the proof. The first claim follows immediately from the
spectral sequence. To prove the second claim, we show that in the
following exact sequence of pro-abelian groups both maps are zero:
$$\{H_q(F(A,I^m))\} \xrightarrow{j}
\{H_q(F(A,I^m),BGL(I^m))\} \xrightarrow{\partial}
\{H_{q-1}(BGL(I^m))\}.$$
We first prove that the map $\partial$ is zero, or equivalently, that
in the following exact sequence of pro-abelian groups, the map $i$ is
injective:
$$\{H_q(F(A,I^m),BGL(I^m))\} \xrightarrow{\partial}
\{H_{q-1}(BGL(I^m))\} \xrightarrow{i}
\{H_{q-1}(F(A,I^m))\}.$$
By hypothesis, the map $i$ is surjective, and therefore, the induced
sequence of $\{ \overline{GL}(A/I^m) \}$-coinvariants again is
exact. Hence, the first statement of the lemma shows that $\partial$
induces an monomorphism of pro-abelian groups
$$\{H_0(\overline{GL}(A/I^m), H_{q-1}(BGL(I^m))) \}
\xrightarrow{\bar{\partial}}
\{H_0(\overline{GL}(A/I^m), H_{q-1}(F(A,I^m))) \}.$$
But, by further hypothesis, the canonical projection
$$\{ H_{q-1}(BGL(I^m)) \} \xrightarrow{\pi'}
\{H_0(\overline{GL}(A/I^m),H_{q-1}(BGL(I^m))) \}$$
is an isomorphism, and hence, the map $\partial$ is injective as
desired.

It remains to show that the map $j$ is zero, or equivalently, that in
the following exact sequence of pro-abelian groups, the $i$ map is
surjective:
$$\{H_q(BGL(I^m))\} \xrightarrow{i} \{H_q(F(A,I^m))\} \xrightarrow{j}
\{H_q(F(A,I^m),BGL(I^m))\}.$$
The map $i$ is surjective by what was proved above, so the induced
sequence of $\{\overline{GL}(A/I^m)\}$-coinvariants again is
exact. Therefore, by the first statement of the lemma, the induced map
$$\{H_0(\overline{GL}(A/I^m), H_q(BGL(I^m)))\} \xrightarrow{\bar{i}}
\{H_0(\overline{GL}(A/I^m), H_q(F(A,I^m)))\}$$
is surjective, so it suffices to show that the canonical projection
$$\{ H_q(F(A,I^m) \} \xrightarrow{\pi}
\{H_0(\overline{GL}(A/I^m),H_q(F(A,I^m)))\}$$
is an isomorphism. But $X(A,I^m)$ is a homotopy pull-back, and hence,
the action of $\overline{GL}(A/I^m) = \pi_1(B\overline{GL}(A/I^m))$ on
$H_q(F(A,I^m))$ is induced from the action of
$\pi_1(B\overline{GL}(A/I^m)^+)$ on $H_q(F(A,I^m))$. Finally, 
$$F(A,I^m) \to BGL(A)^+ \to B\overline{GL}(A/I^m)^+$$
is a fiber sequence of infinite loop spaces, so the fundamental group
of the base acts trivially on the homology of the fiber.
\end{proof}

\begin{prop}\label{homologyexcision}Let $A$ be a unital associative
ring, and let $I \subset A$ be a two-sided ideal. Let $p$ be a prime,
and assume that the pro-abelian group
$$\{ \operatorname{Tor}_q^{\mathbb{Z} \ltimes I^m}
(\mathbb{Z},\mathbb{Z}/p) \}_{m \geqslant 1}$$
is zero, for all positive integers $q$. Then the pro-abelian group
$$\{ H_q(F(A,I^m),BGL(I^m),\mathbb{Z}/p) \}_{m \geqslant 1}$$
is zero, for all integers $q$.
\end{prop}

\begin{proof}We omit the $\mathbb{Z}/p$-coefficients for the
homology. We show that the canonical projection induces an isomorphism
of pro-abelian groups
$$\{H_q(BGL(I^m))\} \xrightarrow{\sim}
\{H_0(\overline{GL}(A/I^m), H_q(BGL(I^m)))\},$$
for all integers $q$. The proposition then follows from
Lemma~\ref{inductionlemma} by induction on $q$ starting from the
trivial case $q = 0$.

It suffices, by the proof of~\cite[Cor.~1.6]{suslinwodzicki}, to show
that for all $m \geqslant 1$, there exists $k \geqslant m$ such that
the image of the map
$$\iota_{k,m*} \colon H_q(GL(I^k)) \to H_q(GL(I^m))$$
is invariant under the adjoint action by
$GL(\mathbb{Z})$. And by~\emph{op.~cit.},~Lemma~1.4, the embedding
$\varphi \colon GL(I^m) \to GL(I^m)$ given by
$$\varphi(\alpha) = \begin{pmatrix}
{1} &
{0} \cr
{0} &
{\alpha} \cr
\end{pmatrix}$$
induces an injection on homology, so it suffices to show that the
image of the composite map
$$H_q(GL(I^k)) \xrightarrow{\iota_{k,m*}} H_q(GL(I^m))
\xrightarrow{\varphi_*} H_q(GL(I^m))$$
is invariant under the adjoint action by the subgroup
$\varphi(GL(\mathbb{Z})) \subset GL(\mathbb{Z})$. In effect, we show
that the image of $\varphi_* \circ \iota_{k,m*}$ is invariant under
the adjoint action by the full group $GL(\mathbb{Z})$. We recall that
$GL(\mathbb{Z})$ is generated by the diagonal matrix
$\operatorname{diag}(-1,1,1,\dots)$ and by the elementary matrices
$e_{s,1}(1)$ and $e_{1,s}(1)$, where $s > 1$. It is clear that
the image of $\varphi_* \circ \iota_{k,m*}$ is invariant under the
action of the diagonal matrix $\operatorname{diag}(-1,1,1,\dots)$. To
show that the image of $\varphi_* \circ \iota_{k,m*}$ is invariant
under the action of the elementary matrices $e_{s,1}(1)$ we consider
the semi-direct product group
$$\Gamma(I^m) = GL(I^m) \ltimes M_{\infty,1}(I^m),$$
where $GL(I^m)$ acts on the group $M_{\infty,1}(I^m)$ of column
vectors by left multiplication. Let $\sigma \colon GL(I^m) \to
\Gamma(I^m)$ and $\pi \colon \Gamma(I^m) \to GL(I^m)$ be the canonical
group homomorphisms given by $\sigma(\alpha) = (\alpha,0)$ and
$\pi(\alpha,v) = \alpha$. Then the map $\varphi$ is equal to the
composition of the group homomorphisms
$$GL(I^m) \xrightarrow{\sigma}
\Gamma(I^m) \xrightarrow{j}
GL(I^m),$$
where
$$j(\alpha,v) = \begin{pmatrix}
{1} &
{0} \cr
{v} &
{\alpha} \cr
\end{pmatrix}.$$
The matrix $e_{s,1}(1) \in GL(\mathbb{Z})$ is the image by $j \colon
\Gamma(\mathbb{Z}) \to GL(\mathbb{Z})$ of a unique matrix
$e_{s,1}'(1) \in \Gamma(\mathbb{Z})$ and
$\pi(e_{s,1}'(1)) \in GL(\mathbb{Z})$ is the identity matrix. It
suffices to show that there exists $k \geqslant m$ such that the image
of the map
$$\iota_{k,m*} \colon H_q(\Gamma(I^k)) \to H_q(\Gamma(I^m))$$
is invariant under the adjoint action of the matrices $e_{s,1}'(1)$,
$s > 1$. It follows from Prop.~\ref{conditionA} below that there
exists $k \geqslant m$ such that
$$\iota_{k,m*} = \iota_{k,m*} \circ \sigma_* \circ \pi_* \colon
H_q(\Gamma(I^k)) \to
H_q(\Gamma(I^m)).$$
Granting this, we obtain that
$$\begin{aligned}
\operatorname{Ad}(e_{s,1}'(1))_* \circ \iota_{k,m*}
{} & = \iota_{k,m*} \circ \operatorname{Ad}(e_{s,1}'(1))_* \cr
{} & = \iota_{k,m*} \circ \sigma_* \circ \pi_* \circ
\operatorname{Ad}(e_{s,1}'(1))_* \cr
{} & = \iota_{k,m*} \circ \sigma_* \circ 
\operatorname{Ad}(\pi(e_{s,1}'(1)))_* \circ \pi_* \cr
{} & = \iota_{k,m*} \circ \sigma_* \circ \pi_* = \iota_{k,m*} \cr
\end{aligned}$$
as desired. An analogous argument shows that there exists (a possibly
larger) $k \geqslant m$ such that the image of $\varphi_* \circ
\iota_{k,m*}$ is invariant under the action of the elementary matrices
$e_{1,s}(1)$.
\end{proof}

\begin{proof}[Proof of Thm.~\ref{ktheory}]It suffices by simple
induction to treat the case $v = 1$. Let $f \colon A \to B$ be as
in the statement. By Prop.~\ref{homologyexcision}, the induced map
$$\{ H_q(F(A,I^m),\mathbb{Z}/p) \}_{m \geqslant 1} \to
\{H_q(F(B,f(I)^m),\mathbb{Z}/p) \}_{m \geqslant 1}$$
is an isomorphism of pro-abelian groups, and~\cite[Cor.~5.8]{panin}
then shows that also the induced map
$$\{ \pi_q(F(A,I^m),\mathbb{Z}/p) \}_{m \geqslant 1} \to
\{ \pi_q(F(B,f(I)^m),\mathbb{Z}/p) \}_{m \geqslant 1}$$
is an isomorphism of pro-abelian groups.
\end{proof}

The remainder of this section is devoted to the proof of
Prop.~\ref{conditionA} below which we used in the proof of
Prop.~\ref{homologyexcision} above. We begin with some general
theory.

Let $A$ be a unital associative ring, and let $I \subset A$ be a
two-sided ideal. Then the powers $\{I^m\}_{m \geqslant 1}$ form a
pro-associative ring without unit. We say that a pro-abelian
group $\{M_m\}_{m \geqslant 1}$ is a \emph{left module} over
$\{I^m\}_{m \geqslant 1}$ if for all $m \geqslant 1$, $M_m$ is a left
$I^m$-module and if for all $k \geqslant m \geqslant 1$, the following
diagram commutes:
$$\xymatrix{
{I^k \otimes M_k} \ar[r]^(.57){\mu_k} \ar[d]^{\iota_{k,m} \otimes
\lambda_{k,m}} &
{M_k} \ar[d]^{\lambda_{k,m}} \cr
{I^m \otimes M_m} \ar[r]^(.57){\mu_m} &
{M_m.} \cr
}$$
By a \emph{homomorphism} from a left $\{I^m\}$-module $\{M_m\}$ to
a left $\{I^m\}$-module $\{M_m'\}$, we mean a strict map of
pro-abelian groups $f \colon \{M_m\} \to \{M_m'\}$ such that for all
$m \geqslant 1$, the map $f_m \colon M_m \to M_m'$ is an $I^m$-module
homomorphism. We say that a left $\{I^m\}$-module $\{P_m\}$ is
\emph{pseudo-free} if for all $m \geqslant 1$, there are $I^m$-module
isomorphisms $\smash{\varphi_m \colon I^m \otimes L_m 
xrightarrow{\sim} P_m}$ with $L_m$ a free abelian group, and if for
all $k \geqslant m \geqslant 1$, there are maps of abelian groups
$\sigma_{k,m} \colon L_k \to L_m$ such that the following diagrams
commute:
$$\xymatrix{
{I^k \otimes L_k} \ar[r]^(.6){\varphi_k} \ar[d]^{\iota \otimes
\sigma_{k,m}} &
{P_k} \ar[d]^{\lambda_{k,m}} \cr
{I^m \otimes L_m} \ar[r]^(.6){\varphi_m} &
{P_m.} \cr
}$$
A homomorphism $f \colon \{P_m\} \to \{M_m\}$ from a pseudo-free left
$\{I^m\}$-module to an arbitrary left $\{I^m\}$-module is said to be

(i) a \emph{special} homomorphism if there exists a strict map of
pro-abelian groups $g \colon \{L_m\} \to \{M_m\}$ such that for all
$m \geqslant 1$, the composite
$$I^m \otimes L_m \xrightarrow{\varphi_m} P_m \xrightarrow{f_m} M_m$$
maps $a \otimes x$ to $a g_m(x)$.

(ii) a \emph{$p^v$-special} homomorphism if the left $\{I^m\}$-module
$\{P_m\}$ admits a strict direct sum decomposition
$$\iota' \oplus \iota'' \colon \{P_m'\} \oplus \{P_m''\} \xrightarrow{\sim}
\{P_m\}$$
such that $\{P_m'\}$ and $\{P_m''\}$ are pseudo-free, such that $f
\circ \iota'$ is special, and such that $f \circ \iota'' = p^v g$, for
some homomorphism $g \colon \{P_m''\} \to \{M_m\}$.

\begin{lem}\label{pseudofreecover}Let $\{M_m\}$ be a left
$\{I^m\}$-module. Then there exists a special homomorphism $f \colon
\{P_m\} \to \{M_m\}$ from a pseudo-free left $\{I^m\}$-module
$\{P_m\}$ such that the image of $f_m$ is equal to $I^m \cdot M_m
\subset M_m$.
\end{lem}

\begin{proof}Let $\mathbb{Z}\{M_m\}$ be the free abelian group
generated by the set of elements of $M_m$. We define $P_m$ to be the
extended left $I^m$-module
$$P_m = I^m \otimes \bigoplus_{s \geqslant m} \mathbb{Z}\{M_s\}$$
and $\lambda_{k,m} \colon P_k \to P_m$ to be the tensor product of the
canonical inclusions. It is clear that $\{P_m\}$ is a pseudo-free left
$\{I^m\}$-module and that the canonical map $f \colon \{P_m\} \to
\{M_m\}$ is a special homomorphism such that the image of $f_m$ is
equal to $I^m \cdot M_m \subset M_m$ as desired.
\end{proof}

Let $\{P_m\}$ be a pseudo-free left $\{I^m\}$-module. By an
\emph{augmented complex} of left $\{I^m\}$-modules over $\{P_m\}$, we
mean a sequence of left $\{I^m\}$-modules and homomorphisms
$$ \dots \xrightarrow{d} \{C_{m,2}\} \xrightarrow{d}
\{C_{m,1}\} \xrightarrow{d} \{C_{m,0}\} \xrightarrow{\epsilon}
\{P_m\}$$
such that $d_m \circ d_m$ and $\epsilon_m \circ d_m$ are equal to
zero, for all $m \geqslant 1$. We say that the augmented complex
$\epsilon \colon \{C_{m,*}\} \to \{P_m\}$ is

(i) \emph{special} if the $\{C_{m,q}\}$ are pseudo-free, if $\epsilon$
is special, and if the homology of the augmented complex $\epsilon_m
\colon C_{m,*} \to P_m$ is annihilated by $I^m$.

(ii) \emph{$p^v$-special} if the $\{C_{m,q}\}$ are pseudo-free, if
$\epsilon$ is $p^v$-special, and if the homology of the augmented
complex $\epsilon_m \colon C_{m,*} \to P_m$ is annihilated both by
$I^m$ and by $p^v$.

We remark that, for a pseudo-free left $\{I^m\}$-module $\{P_m\}$, the
identity homomorphism is usually not special, and hence, does not
constitute an augmented special complex.

\begin{cor}\label{specialcomplex}Let $\{P_m\}$ be a pseudo-free left
$\{I^m\}$-module. Then there exists both a special augmented complex
and a $p^v$-special augmented complex over $\{P_m\}$.
\end{cor}

\begin{proof}It is proved by an easy induction argument based on
Lemma~\ref{pseudofreecover} that there exists a special augmented
complex
$$\xymatrix{
{ \dots } \ar[r]^{d} &
{ \{C_{m,2}\} } \ar[r]^{d} &
{ \{C_{m,1}\} } \ar[r]^{d} &
{ \{C_{m,0}\} } \ar[r]^{\epsilon} &
{ \{P_m\}. } \cr
}$$
Then the total complex of the double-complex
$$\xymatrix{
{ \dots } \ar[r]^{d}  &
{ \{C_{m,2}\} } \ar[r]^{d} \ar[d]^{p^v} &
{ \{C_{m,1}\} } \ar[r]^{d} \ar[d]^{p^v} &
{ \{C_{m,0}\} } \ar[r]^{\epsilon} \ar[d]^{p^v} &
{ \{P_m\} } \ar[d]^{p^v} \cr
{ \dots } \ar[r]^{d} &
{ \{C_{m,2}\} } \ar[r]^{d} &
{ \{C_{m,1}\} } \ar[r]^{d} &
{ \{C_{m,0}\} } \ar[r]^{\epsilon} &
{ \{P_m\} } \cr
}$$
is a $p^v$-special augmented complex.
\end{proof}

\begin{prop}\label{specialresolution}Let $A$ be a unital associative
ring, and let $I \subset A$ be a two-sided ideal. Let $p$ be a prime, 
and assume that for every positive integer $q$, the following
pro-abelian group is zero:
$$\{\operatorname{Tor}_q^{\mathbb{Z} \ltimes I^m}
(\mathbb{Z},\mathbb{Z}/p)\}_{m \geqslant 1}.$$
Then every $p^v$-special augmented complex of $\{I^m\}$-modules
$$ \dots \xrightarrow{d} \{C_{m,2}\} \xrightarrow{d} \{C_{m,1}\}
\xrightarrow{d} \{C_{m,0}\} \xrightarrow{\epsilon} \{P_m\},$$
is exact as an augmented complex of pro-abelian groups.
\end{prop}

\begin{proof}It follows by simple induction that the pro-abelian group
$$\{ \operatorname{Tor}_q^{\mathbb{Z} \ltimes I^m}
(\mathbb{Z},\mathbb{Z}/p^v) \}_{m \geqslant 1}$$
is zero, for all positive integers $q$ and $v$. We show by induction
on the degree that the following complex of pro-abelian groups has
zero homology:
$$ \dots \xrightarrow{d} \{C_{m,2}\} \xrightarrow{d} \{C_{m,1}\}
\xrightarrow{d} \{C_{m,0}\} \xrightarrow{d} \{C_{m,-1}\} = \{P_m\}.$$
We let $Z_{m,q}$ and $B_{m,q-1}$ be the kernel and image of $d_m
\colon C_{m,q} \to C_{m,q-1}$. Since $I^m$ and $p^v$ both annihilate
$H_q(C_{m,*}) = Z_{m,q}/B_{m,q}$, we have
$$I^m Z_{m,q} + p^vZ_{m,q} \subset B_{m,q}.$$
We first consider $q = -1$. Then $Z_{m,-1} = P_m$, and hence
$I^mP_m + p^vP_m$ is contained in $B_{m,-1}$. But the definition of a
$p^v$-special homomorphism shows that also $B_{m,-1} \subset I^mP_m +
p^vP_m$. It follows that we have an isomorphism
$$H_{-1}(C_{m,*}) \xleftarrow{\sim}
\operatorname{Tor}_0^{\mathbb{Z} \ltimes P_m}(P_m,\mathbb{Z}/p^v).$$
Now, since $\{P_m\}$ is a pseudo-free left $\{I^m\}$-module, we
have an isomorphism of left $I^m$-modules $I^m \otimes L_m
\xrightarrow{\sim} P_m$, with $L_m$ a free abelian group, so
$$\begin{aligned}
H_{-1}(C_{m,*}) & \xleftarrow{\sim}
\operatorname{Tor}_0^{\mathbb{Z} \ltimes I^m}
(P_m,\mathbb{Z}/p^v) \xleftarrow{\sim}
\operatorname{Tor}_0^{\mathbb{Z} \ltimes I^m}
(I^m \otimes L_m,\mathbb{Z}/p^v) \cr
{} &  \xleftarrow{\sim} \operatorname{Tor}_0^{\mathbb{Z} \ltimes I^m}
(I^m,\mathbb{Z}/p^v) \otimes L_m \xleftarrow{\sim}
\operatorname{Tor}_1^{\mathbb{Z} \ltimes I^m}
(\mathbb{Z},\mathbb{Z}/p^v) \otimes L_m. \cr
\end{aligned}$$
This show that the pro-abelian group $\{H_{-1}(C_{m,*})\}$ is zero as
stated.

Next, we assume, inductively, that the pro-abelian group
$\{H_i(C_{m,*})\}$ is zero, for $i < q$, and show that the pro-abelian
group $\{H_q(C_{m,*})\}$ is zero. The inductive hypothesis is
equivalent to the statement that the following sequence of pro-abelian
groups is exact:
$$0 \to \{Z_{m,q}\} \to \{C_{m,q}\} \to \dots \to \{C_{m,0}\} \to
\{P_m\} \to 0.$$
This gives rise to the following spectral sequence of pro-abelian
groups:
$$\begin{aligned}
E_{s,t}^1 & = \begin{cases}
\{\operatorname{Tor}_t^{\mathbb{Z} \ltimes I^m}
(C_{m,s},\mathbb{Z}/p^v)\}, &
\text{if $0 \leqslant s \leqslant q$,} \cr
\{\operatorname{Tor}_t^{\mathbb{Z} \ltimes I^m}
(Z_{m,q},\mathbb{Z}/p^v)\}, &
\text{if $s = q + 1$,} \cr
\{0\}, & \text{if $s > q + 1$,} \cr
\end{cases} \cr
{} & \Rightarrow \{\operatorname{Tor}_{s+t}^{\mathbb{Z} \ltimes I^m}
(P_m,\mathbb{Z}/p^v) \}. \cr
\end{aligned}$$
Since the left $\{I^m\}$-modules $\{C_{m,s}\}$ are pseudo-free, we see
as before that, for $0 \leqslant s \leqslant q$, the pro-abelian
groups $E_{s,t}^1$ are zero. It follows that the edge-homomorphism of
the spectral sequence
$$\{\operatorname{Tor}_{q+1}^{\mathbb{Z} \ltimes I^m}
(P_m,\mathbb{Z}/p^v)\} \to
\{\operatorname{Tor}_0^{\mathbb{Z} \ltimes I^m}
(Z_{m,q},\mathbb{Z}/p^v)\}$$
is an isomorphism of pro-abelian groups. The pro-abelian group on the
left-hand side is zero, since the left $\{I^m\}$-module $\{P_m\}$ is
pseudo-free. Hence, we may conclude that the following pro-abelian
group is zero:
$$\{\operatorname{Tor}_0^{\mathbb{Z} \ltimes I^m}
(Z_{m,q},\mathbb{Z}/p^v)\} = \{Z_{m,q}/(I^m Z_{m,q} + p^vZ_{m,q}) \}.$$
Finally, since $I^m Z_{m,q} + p^vZ_{m,q}$ is contained in $B_{m,q}$,
this pro-abelian group surjects onto the pro-abelian group $\{Z_{m,q}
/ B_{m,q}\} = \{H_q(C_{m,*})\}$ which therefore is zero as stated.
\end{proof}

\begin{prop}\label{pseudofreezero}Let $A$ be a unital associative
ring and let $I \subset A$ be a two-sided ideal. Let $p$ be a prime
and assume that for every positive integer $q$, the following
pro-abelian group is zero:
$$\{\operatorname{Tor}_q^{\mathbb{Z} \ltimes I^m}
(\mathbb{Z},\mathbb{Z}/p)\}_{m \geqslant 1}.$$
Suppose that $F_i$, $i \geqslant 0$, is a sequence of possibly
non-additive functors from the category of $\{I^m\}$-modules to the
category of pro-abelian groups that satisfy the following conditions
(i)---(iii):

(i) For every $\{I^m\}$-module $\{M_m\}$, the pro-abelian group
$F_0(\{M_m\})$ is zero.

(ii) If $\epsilon \colon \{M_m[-]\} \to \{M_m\}$ is an augmented
simplicial $\{I^m\}$-module such that the associated augmented
complex of pro-abelian groups is exact, then there is a spectral
sequence of pro-abelian groups
$$E_{s,t}^1 = F_t(\{M_m[s]\}) \Rightarrow F_{s+t}(\{M_m\})$$
such that the edge-homomorphism $F_t(\{M_m[0]\}) \to F_t(\{M_m\})$ is
the map induced by the augmentation $\epsilon \colon \{M_m[0]\} \to
\{M_m\}$.

(iii) For every $i \geqslant 0$, there exists $v \geqslant 1$ such
that for every $p^v$-special homomorphism $f \colon \{P_m\} \to
\{M_m\}$, the induced map of pro-abelian groups
$$f_* \colon F_i(\{P_m\}) \to F_i(\{M_m\})$$
is zero.

Then for every pseudo-free $\{I^m\}$-module $\{P_m\}$ and for every
$i \geqslant 0$, the pro-abelian group $F_i(\{P_m\})$ is zero.
\end{prop}

\begin{proof}The proof is by induction on $i \geqslant 0$. The case
$i = 0$ is trivial, so we assume the statement for $i-1$ and prove it
for $i$. We choose $v$ such that for every $0 \leqslant t \leqslant i$
and every $p^v$-special homomorphism $f \colon \{P_m\} \to \{M_m\}$,
the induced map of pro-abelian groups
$$f_* \colon F_t(\{P_m\}) \to F_t(\{M_m\})$$
is zero. We let $\{P_m\}$ be a pseudo-free $\{I^m\}$-module and
choose a $p^v$-special augmented complex of $\{I^m\}$-modules
$\epsilon \colon \{C_{m,*}\} \to \{P_m\}$. It follows from
Prop.~\ref{specialresolution} that this augmented complex is exact as
an augmented complex of pro-abelian groups. We let $\epsilon \colon
\{P_m[-]\} \to \{P_m\}$ be the augmented simplicial $\{I^m\}$-module
associated with the augmented complex of $\{I^m\}$-modules $\epsilon
\colon \{C_{m,*}\} \to \{P_m\}$ by the Dold-Kan
correspondence~\cite{weibel1}. By assumption, we then have the 
following spectral sequence of pro-abelian groups:
$$E_{s,t}^1 = F_t(\{P_m[s]\}) \Rightarrow F_{s+t}(\{P_m\}).$$
The inductive hypothesis shows that the pro-abelian groups $E_{s,t}^1$
are zero, for all $s$ and all $0 \leqslant t < i$. Hence, the
edge-homomorphism 
$$F_i(\{P_m[0]\}) \to F_i(\{P_m\})$$
is an epimorphism of pro-abelian groups. But the edge-homomorphism is
equal to the map induced by the special homomorphism $\epsilon \colon
\{P_m[0]\} \to \{P_m\}$ and therefore is zero by assumption. This
completes the proof.
\end{proof}

If a group $G$ acts from the left on a group $H$, the semi-direct
product group $G \ltimes H$ is defined to be the set $G \times H$ with
multiplication
$$(g,h) \cdot (g',h') = (gg', (g'{}^{-1}h)h').$$
The projection $\pi \colon G \ltimes H \to G$ and the section $\sigma
\colon G \to G \ltimes H$ defined by $\sigma(g) = (g,1)$ are both
group-homomorphisms, and hence, the homology of the semi-direct
product group decomposes canonically as a direct sum
$$H_q(G \ltimes H,\mathbb{Z}/p) \xleftarrow{\sim}
H_q(G,\mathbb{Z}/p) \oplus H_q(G \ltimes H, G, \mathbb{Z}/p).$$
We write $\mathbb{Z}[X]$ for the free abelian group generated by the
set $X$.

\begin{prop}\label{spectralsequence}Let $\{G_m\}$ be a pro-group, 
let $\epsilon \colon \{H_m[-]\} \to \{H_m\}$ be a strict augmented
simplicial pro-group with a strict left $\{G_m\}$-action, and suppose
that the associated complex of pro-abelian groups
$$\dots \to \{\mathbb{Z}[H_m[2]]\} \to \{\mathbb{Z}[H_m[1]]\} \to
\{\mathbb{Z}[H_m[0]]\} \xrightarrow{\epsilon}
\{\mathbb{Z}[H_m]\} \to 0$$
is exact. Then there is a natural spectral sequence of pro-abelian
groups
$$E_{s,t}^1 = \{H_t(G_m \ltimes H_m[s], G_m, \mathbb{Z}/p) \}
\Rightarrow \{H_{s+t}(G_m \ltimes H_m, G_m, \mathbb{Z}/p) \},$$
and the edge-homomorphism
$$\{H_t(G_m \ltimes H_m[0], G_m, \mathbb{Z}/p)\} \to
\{H_t(G_m \ltimes H_m, G_m, \mathbb{Z}/p)\}$$
is equal to the map induced from $\epsilon$.
\end{prop}

\begin{proof}Suppose first that $G$ is a group and that $M_*$ is a
chain complex of left $G$-modules. The hyper-homology of $G$ with
coefficients in $M_*$ is defined to be the homology $H_*(G;M_*)$ of
the  total complex $C_*(G;M_*)$ of the double-complex defined by the
bar-construction $B_*(\mathbb{Z}[G],\mathbb{Z}[G],M_*)$. The
filtration in the bar-direction gives rise to a functorial spectral
sequence
$$E_{s,t}^2 = H_s(G; H_t(M_*)) \Rightarrow H_{s+t}(G;M_*),$$
where $H_t(M_*)$ is considered as a complex concentrated in degree
zero.

If $K$ is a subgroup of $G$, the complex $C_*(G,K;M_*)$ is defined to
be the mapping cone of the map $C_*(K;M_*) \to C_*(G;M_*)$ induced by
the inclusion. We write $\tilde{C}_*(G;M_*)$ instead of
$C_*(G,\{1\};M_*)$. Suppose that $G$ acts from the left on a group
$H$, and that $M_*$ is a complex of left $G \ltimes H$-modules. Then
there is a canonical isomorphism of complexes
$$C_*(G \ltimes H; M_*) \xleftarrow{\sim} C_*(G; C_*(H;M))$$
where the $G$-action on $C_*(H; M_*)$ is induced from the actions by
$G$ on $H$ and $M_*$. This isomorphism induces an isomorphism of
complexes
$$C_*(G \ltimes H, G; M_*) \xleftarrow{\sim}
C_*(G; \tilde{C}_*(H; M_*)).$$
In the following, we consider this isomorphism only for $M_* =
\mathbb{Z}/p$.

Let $\epsilon \colon H[-] \to H$ be an augmented simplicial group with
a left action by the group $G$, and let
$\smash{\tilde{C}_*(H[-];\mathbb{Z}/p)_*}$ be the total-complex of the
simplicial complex $\tilde{C}_*(H[-];\mathbb{Z}/p)$. Then the
filtration in the simplicial direction gives rise to a functorial
spectral sequence of hyper-homology groups
$$E_{s,t}^1 = H_t(G; \tilde{C}_*(H[s];\mathbb{Z}/p)) \Rightarrow
H_{s+t}(G; \tilde{C}_*(H[-];\mathbb{Z}/p)_*).$$
Suppose now that $\{G_m\}$ and $\epsilon \colon \{H_m[-]\} \to
\{H_m\}$ are as in the statement. Then we obtain a spectral sequence
of strict pro-abelian groups
$$E_{s,t}^1 = \{H_t(G_m; \tilde{C}_*(H_m[s];\mathbb{Z}/p))\}
\Rightarrow \{H_{s+t}(G_m; \tilde{C}_*(H_m[-];\mathbb{Z}/p)_*)\}$$
which we claim has the desired form. Indeed, the isomorphism of
complexes that we established above identifies
$$\{H_t(G_m; \tilde{C}_*(H_m[s];\mathbb{Z}/p))\} \xrightarrow{\sim}
\{H_t(G_m \ltimes H_m[s],G_m;\mathbb{Z}/p)\}$$
and
$$\{H_{s+t}(G_m; \tilde{C}_*(H_m;\mathbb{Z}/p))\} \xrightarrow{\sim}
\{H_{s+t}(G_m \ltimes H_m,G_m;\mathbb{Z}/p)\}.$$
Hence, it remains only to show that the map induced by the
augmentation
$$\epsilon_* \colon \{H_{s+t}(G_m; \tilde{C}_*
(H_m[-];\mathbb{Z}/p)_*)\} \to
\{H_{s+t}(G_m; \tilde{C}_*(H_m;\mathbb{Z}/p)) \}$$
is an isomorphism of pro-abelian groups. But the hypothesis of the
statement implies that for all $t$, the map induced by the
augmentation
$$\epsilon_* \colon \{H_t(\tilde{C}_*(H_m[-];\mathbb{Z}/p)_*) \} \to
\{\tilde{H}_t(H_m;\mathbb{Z}/p)\}$$
is an isomorphism of pro-abelian groups, and the desired isomorphism
then follows from the spectral sequence at the beginning of the proof.
\end{proof}

Let again $A$ be a unital associative ring, let $I \subset A$ be a
two-sided ideal, and let $\{M_m\}$ be a left $\{I^m\}$-module. Then
the group of column vectors $M_{\infty,1}(M_m)$ is naturally a left
$GL(I^m)$-module, and we consider the canonical inclusion of $GL(I^m)$
in the semi-direct product $GL(I^m) \ltimes M_{\infty,1}(M_m)$. The
relative homology groups form a sequence of non-additive functors
$$F_i(\{M_m\}) = \{H_i(GL(I^m) \ltimes M_{\infty,1}(M_m),
GL(I^m),\mathbb{Z}/p)\}$$
from the category of left $\{I^m\}$-modules to the category of
pro-abelian groups. This sequence of functors trivially satisfies the
condition~(i) of Prop.~\ref{pseudofreezero} and it follows from
Prop.~\ref{spectralsequence} that it also satisfies the
condition~(ii). Finally, it follows from~\cite[Cor.~6.6]{suslin4} that
it satisfies the condition~(iii). We now prove the following result
which we used in the proof of Prop.~\ref{homologyexcision}:

\begin{prop}\label{conditionA}Let $A$ be a unital associative ring and
let $I \subset A$ be a two-sided ideal. Let $p$ be a prime and assume
that for every positive integer $q$, the following pro-abelian group
is zero:
$$\{\operatorname{Tor}_q^{\mathbb{Z} \ltimes I^m}
(\mathbb{Z},\mathbb{Z}/p)\}_{m \geqslant 1}.$$
Then for all integers $q$, the pro-abelian group
$$\{H_q(GL(I^m) \ltimes M_{\infty,1}(I^m),GL(I^m),\mathbb{Z}/p)\}$$
is zero.
\end{prop}

\begin{proof}Indeed, this follows from Prop.~\ref{pseudofreezero},
since the $\{I^m\}$-module $\{I^m\}$ is pseudo-free.
\end{proof}

\section{Topological cyclic homology}\label{thhsection}

This section is devoted to the proof of the following analog of
Thm~\ref{ktheory} for topological cyclic homology.

\begin{thm}\label{maintc}Let $f \colon A \to B$ be a map of
unital associative rings, and let $I \subset A$ be a two-sided ideal
such that $f \colon I \to f(I)$ is an isomorphism onto a two-sided
ideal of $B$. Let $p$ be a prime and suppose that for all positive
integers $q$, the pro-abelian group
$$\{ \operatorname{Tor}_q^{\mathbb{Z} \ltimes I^m}
(\mathbb{Z},\mathbb{Z}/p) \}_{m \geqslant 1}$$
is zero. Then for all integers $q$, and all positive integers $n$ and
$v$, the induced map of pro-abelian groups
$$\{ \operatorname{TC}_q^n(A,I^m;p,\mathbb{Z}/p^v) \}_{m \geqslant 1} \to
\{ \operatorname{TC}_q^n(B,f(I)^m;p,\mathbb{Z}/p^v) \}_{m \geqslant 1}$$
is an isomorphism.
\end{thm}

The proof is similar, in outline, to the original proof by
Wodzicki~\cite{wodzicki} of excision in cyclic homology for $H$-unital
rings, but the technical details are somewhat different. We first
show, as in the original case, that the statement for topological
cyclic homology follows from the analogous statement for topological
Hochschild homology. We briefly recall topological Hochschild homology
and topological cyclic homology. The reader is referred
to~\cite[Sect.~1]{hm4} and~\cite{h3} for a more detailed account.

The topological Hochschild spectrum $T(A)$ of a unital associated ring
$A$ is a \emph{cyclotomic spectrum} in the sense
of~\cite[Def.~2.2]{hm}. In particular, $T(A)$ is an object of the
$\mathbb{T}$-stable homotopy category, where $\mathbb{T}$ is the
circle group, and is equipped with additional cyclotomic structure
maps which on the level of the equivariant homotopy groups
$$\operatorname{TR}_q^n(A;p) = [S^q \wedge (\mathbb{T}/C_{p^{n-1}})_+,
T(A)]_{\mathbb{T}}$$
give rise to the so-called~\emph{restriction} maps
$$R \colon \operatorname{TR}_q^n(A;p) \to
\operatorname{TR}_q^{n-1}(A;p).$$
We consider the family of groups
$\smash{\{\operatorname{TR}_q^n(A;p)\}_{n \geqslant 1}}$ as a
pro-abelian group with the restriction maps as the structure maps. As
a consequence of the cyclotomic spectrum structure, there are natural
long-exact sequences
$$\cdots \to \mathbb{H}_q(C_{p^{n-1}}, T(A)) \to
\operatorname{TR}_q^n(A;p) \xrightarrow{R}
\operatorname{TR}_q^{n-1}(A;p) \to \cdots.$$
The group $\operatorname{TR}_q^1(A;p)$ is the $q$th topological
Hochschild group $\operatorname{TH}_q(A)$, and there are and natural
spectral sequences
$$E_{s,t}^2 = H_s(C_{p^{n-1}}, \operatorname{TH}_t(A)) \Rightarrow 
\mathbb{H}_{s+t}(C_{p^{n-1}}, T(A)).$$
The canonical projection $\mathbb{T}/C_{p^{n-2}} \to
\mathbb{T}/C_{p^{n-1}}$ induces the~\emph{Frobenius} map
$$F \colon \operatorname{TR}_q^n(A;p) \to
\operatorname{TR}_q^{n-1}(A;p),$$
and the topological cyclic homology groups
$\operatorname{TC}_q^n(A;p)$ are defined so that there is a natural
long-exact sequence of abelian groups
$$\cdots \to \operatorname{TC}_q^n(A;p) \to
\operatorname{TR}_q^n(A;p) \xrightarrow{R-F}
\operatorname{TR}_q^{n-1}(A;p) \xrightarrow{\partial}
\operatorname{TC}_{q-1}^n(A;p) \to \cdots.$$
Let $I \subset A$ be a two-sided ideal. Then one has the following
natural distinguished triangle in the $\mathbb{T}$-stable category:
$$T(A,I) \to T(A) \to T(A/I) \xrightarrow{\partial} \Sigma T(A,I).$$
Moreover, B\"okstedt and Madsen~\cite[Prop.~10.2]{bokstedtmadsen}
have given an explicit construction of this distinguished triangle
with the property that it is also a sequence of cyclotomic spectra. It
follows that there are induced natural long-exact sequences
$$\cdots \to \operatorname{TC}_q^n(A,I;p) \to
\operatorname{TC}_q^n(A;p) \to
\operatorname{TC}_q^n(A/I;p) \xrightarrow{\partial}
\operatorname{TC}_{q-1}^n(A,I;p) \to \cdots$$
and that there are relative versions of the long-exact sequences and
the spectral sequence considered above. Hence, to prove
Thm.~\ref{maintc} above, it suffices to prove the following analogous
result for topological Hochschild homology. We remark that this
statement concerns only the underlying non-equivariant spectrum
$\operatorname{TH}(A,I)$ of the $\mathbb{T}$-equivariant spectrum
$T(A,I)$.

\begin{thm}\label{mainthh}Let $f \colon A \to B$ be a map of
unital associative rings, and let $I \subset A$ be a two-sided ideal
such that $f \colon I \to f(I)$ is an isomorphism onto a two-sided
ideal of $B$. Let $p$ be a prime, and suppose that for all positive
integers $q$, the pro-abelian group
$$\{ \operatorname{Tor}_q^{\mathbb{Z} \ltimes I^m}
(\mathbb{Z},\mathbb{Z}/p) \}_{m \geqslant 1}$$
is zero. Then for all integers $q$, the induced map
$$\{ \operatorname{TH}_q(A,I^m,\mathbb{Z}/p) \}_{m \geqslant 1} \to
\{ \operatorname{TH}_q(B,f(I)^m,\mathbb{Z}/p) \}_{m \geqslant 1}$$
is an isomorphism of pro-abelian groups.
\end{thm}

The proof of Thm.~\ref{mainthh} is rather lengthy and a brief
outline is in order. The definition of topological Hochschild homology
extends to non-unital rings and there is a canonical map
$\operatorname{TH}(I) \to \operatorname{TH}(A,I)$. We show that the
mapping cone of this map is
$B(I)$-cellular~\cite[Def.~5.1.4]{hirschhorn}, where $B(I)$ is the 
bar-construction of the Eilenberg-MacLane spectrum of $I$. The precise
statement is given in Lemma~\ref{gradedpieces} below. Finally, we
relate the $\mathbb{Z}/p$-homology groups of the spectrum $B(I)$ to
the groups $\smash{\operatorname{Tor}_q^{\mathbb{Z} \ltimes I}
(\mathbb{Z},\mathbb{Z}/p)}$ by means of a spectral sequence. We then
conclude that for all integers $q$, the canonical map
$$\{ \operatorname{TH}_q(I^m,\mathbb{Z}/p) \}_{m \geqslant 1} \to
\{ \operatorname{TH}_q(A,I^m,\mathbb{Z}/p) \}_{m \geqslant 1}$$
is an isomorphism of pro-abelian groups. Hence, the right-hand
side is independent of the ring $A$ as stated in Thm.~\ref{mainthh}.

We recall from~\cite{hoveyshipleysmith} that the category of symmetric
spectra is a closed model category and that the associated homotopy
category is the stable homotopy category. Moreover, the smash product
of symmetric spectra defines a closed symmetric monoidal structure
which induces a closed symmetric monoidal structure on the stable
homotopy category.

We will consider various symmetric spectra that are defined as the
realization of semi-simplicial symmetric spectra. To briefly recall
this, let $\mathbf{\Delta}_{\text{mon}}$ be the category whose
objects are the finite ordered sets
$$[k] = \{1,2,\dots,k\}, \hskip4mm k \geqslant 0,$$
and whose morphisms are all \emph{injective} order-preserving maps. 
The morphisms are generated by the coface maps $d^r \colon [k-1] \to
[k]$, $0 \leqslant r \leqslant k$, defined by $d^r(i) = i$, for $0
\leqslant i < r$, and $d^r(i) = i+1$, for $r \leqslant i < k$. Then a
semi-simplicial symmetric spectrum is a functor $X[-]$ from
$\mathbf{\Delta}_{\text{mon}}^{\operatorname{op}}$ to the category of
symmetric spectra, and the realization is the symmetric spectrum given
by the coequalizer
\begin{equation}\label{realization}
\xymatrix{
{\displaystyle{\bigvee X[k'] \wedge (\Delta^k)_+}}
\ar[r]<.8ex> \ar[r]<-.8ex> &
{\displaystyle{\bigvee X[k] \wedge (\Delta^k)_+}} \ar[r] &
{\big{\|} [k] \mapsto X[k] \big{\|},} \cr
}
\end{equation}
where the wedge-sums in the middle and on the left-hand side range
over all objects and morphisms, respectively, of
$\mathbf{\Delta}_{\text{mon}}$. The two maps are given, on the summand
indexed by the morphism $\theta \colon [k] \to
[k']$, by the maps $X[k'] \wedge \theta_*$ and $\theta^* \wedge
(\Delta^k)_+$, respectively. The skeleton filtration of the
realization gives rise to a natural spectral sequence
\begin{equation}\label{skeletonspectralsequence}
E_{s,t}^1 = \pi_t(X[s]) \Rightarrow \pi_{s+t}(\big{\|} [k] \mapsto X[k]
\big{\|})
\end{equation}
where the $d^1$-differential is given by the alternating sum of the
maps induced by the face maps $d_r \colon X[s] \to X[s-1]$,
$0 \leqslant r \leqslant s$.

Let $I$ be a non-unital associative ring. We abuse notation and write
$I$ also for the associated Eilenberg-MacLane
spectrum~\cite[Ex.~1.2.5]{hoveyshipleysmith}. The ring structure gives
rise to a map of symmetric spectra
$$\mu \colon I \wedge I \to I$$
that makes $I$ a non-unital associative symmetric ring spectrum. We
define the topological Hochschild spectrum of $I$ to be the symmetric
spectrum
$$\operatorname{TH}(I) = \big{\|} [k] \mapsto
\operatorname{TH}(I)[k] \big{\|}$$
that is the realization of the semi-simplicial symmetric spectrum
defined as follows. The symmetric spectrum in degree $k$ is the
$(k+1)$-fold smash product
$$\operatorname{TH}(I)[k] = I^{\wedge (k+1)}$$
and the face maps $d_r \colon \operatorname{TH}(I)[k] \to
\operatorname{TH}(I)[k-1]$, $0 \leqslant r \leqslant k$, are defined
by
$$d_r = \begin{cases}
I^{\wedge r} \wedge \mu \wedge I^{\wedge (k-1-r)} &
0 \leqslant r < k, \cr
(\mu \wedge I^{\wedge (k-1)}) \circ \tau & r = k. \cr
\end{cases}$$
Here $\tau$ is the canonical isomorphism that switches the first $k$
smash factors and the last smash factor. We recall
from~\cite[Thm.~4.2.8]{shipley} that if $A$ is a unital associative
ring, the symmetric spectrum $\operatorname{TH}(A)$ defined here is
canonically weakly equivalent to the topological Hochschild spectrum
defined by B\"okstedt. We note, however, that the symmetric spectrum
$\operatorname{TH}(A)$ defined here is \emph{not} equipped with the
additional cyclotomic structure that is needed for the definition of
topological cyclic homology. But this extra structure is not relevant
for the proof of Thm.~\ref{mainthh}.

Let $A$ be a unital associative ring, and let $I \subset A$ be a
two-sided ideal. We continue to write $A$ and $I$ for the associated
Eilenberg-MacLane spectra. Then $A$ is a unital associative symmetric
ring spectrum, $I$ is a unital and symmetric $A$-$A$-bimodule spectrum,
and the map $\rho \colon I \to A$ induced by the inclusion of $I$ in
$A$ is a map of $A$-$A$-bimodule spectra. We follow
B\"{o}kstedt-Madsen~\cite[Appendix]{bokstedtmadsen} and define a
symmetric spectrum
$$\operatorname{TH}(A)' =
\big{\|} [k] \mapsto \operatorname{TH}(A)'[k] \big{\|}$$ 
as the realization of a semi-simplicial symmetric spectrum that we now
define. Let $\mathcal{P}[k]$ be the category with objects all tuples
$\epsilon = (\epsilon_0,\dots,\epsilon_k)$ of integers $\epsilon_i \in
\{0,1\}$ and with a single morphism from $\epsilon$ to $\epsilon'$ if
$\epsilon_i \geqslant \epsilon_i'$, for all $0 \leqslant i \leqslant
k$. So $\mathcal{P}[k]$ is an index category for $(k+1)$-dimensional
cubical diagrams. We let $E[k]$ be the functor from $\mathcal{P}[k]$ to
the category of symmetric spectra that is given on objects and
morphisms, respectively, by
$$\begin{aligned}
E[k](\epsilon_0,\dots,\epsilon_r,\dots,\epsilon_k) & =
E^{\epsilon_0} \wedge \dots \wedge E^{\epsilon_r} \wedge \dots
\wedge E^{\epsilon_k} \cr
E[k](\epsilon_0,\dots,\iota_r,\dots,\epsilon_k) & = 
E^{\epsilon_0} \wedge \dots \wedge \rho \wedge \dots \wedge
E^{\epsilon_k}, \cr
\end{aligned}$$
where $E^0 = A$ and $E^1 = I$, and where
$(\epsilon_0,\dots,\iota_r,\dots,\epsilon_k)$ denotes the unique
morphism from $(\epsilon_0,\dots,1,\dots,\epsilon_k)$ to
$(\epsilon_0,\dots,0,\dots,\epsilon_k)$. Then
$$\operatorname{TH}(A)'[k] =
\operatornamewithlimits{hocolim}_{\mathcal{P}[k]} E[k];$$
see~\cite[Chap.~XII]{bousfieldkan} for the definition of homotopy
colimits. To define the face maps, we first note that the
ring and bimodule spectrum maps define maps of symmetric spectra
$\mu_{\epsilon,\epsilon'} \colon E^{\epsilon} \wedge E^{\epsilon'} \to
E^{\epsilon \sqcup \epsilon'}$ where $\epsilon
\sqcup \epsilon' = \max \{ \epsilon, \epsilon' \}$. Then the face maps
$$d_r \colon \operatorname{TH}(A)'[k] \to \operatorname{TH}(A)'[k-1],
\hskip3mm 0 \leqslant r \leqslant k,$$
are defined to be the composites
$$\operatornamewithlimits{hocolim}_{\mathcal{P}[k]} E[k]
\xrightarrow{\delta_r} 
\operatornamewithlimits{hocolim}_{\mathcal{P}[k]} E[k-1] \circ \mathfrak{d}_r 
\to
\operatornamewithlimits{hocolim}_{\mathcal{P}[k-1]} E[k-1],$$
where $\mathfrak{d}_r \colon \mathcal{P}[k] \to \mathcal{P}[k-1]$ is
the functor
$$\mathfrak{d}_r(\epsilon_0,\dots,\epsilon_k) = \begin{cases}
(\epsilon_0,\dots, \epsilon_r \sqcup \epsilon_{r+1}, \epsilon_k) & 0
\leqslant r < k, \cr
(\epsilon_k \sqcup \epsilon_0, \epsilon_1, \dots, \epsilon_{k-1}) &
r=k, \cr
\end{cases}$$
and where $\delta_r \colon E[k] \to E[k-1] \circ \mathfrak{d}_r$ is the
natural transformation
$$\delta_r = \begin{cases}
E^{\epsilon_0} \wedge \dots \wedge E^{\epsilon_{r-1}} \wedge
\mu_{\epsilon_r,\epsilon_{r+1}} \wedge E^{\epsilon_{r+2}} \wedge \dots
\wedge E^{\epsilon_k} & 0 \leqslant r < k, \cr
(\mu_{\epsilon_k,\epsilon_0} \wedge E^{\epsilon_1} \wedge \dots \wedge
E^{\epsilon_{k-1}}) \circ \tau & r = k. \cr
\end{cases}$$
The object $(0,\dots,0)$ in $\mathcal{P}[k]$ is a final object, and
hence, the canonical map
$$\operatorname{TH}(A)[k] = E[k](0,\dots,0) \to
\operatornamewithlimits{hocolim}_{\mathcal{P}[k]} E[k] =
\operatorname{TH}(A)'[k]$$
is a weak equivalence~\cite[XII.3.1]{bousfieldkan}. The canonical maps
form a map of semi-simplicial symmetric spectra and the spectral
sequence~(\ref{skeletonspectralsequence}) shows that the induced map
of realizations
$$\operatorname{TH}(A) \to \operatorname{TH}(A)'$$
is a weak equivalence of symmetric spectra.

Let $\mathcal{P}'[k] \subset \mathcal{P}[k]$ be the full
subcategory whose objects are the tuples $\epsilon$ such that not all
of $\epsilon_0,\dots,\epsilon_k$ are equal to zero. The symmetric
spectrum
$$\operatorname{TH}(A,I)[k] =
\operatornamewithlimits{hocolim}_{\mathcal{P}'[k]} E[k]$$
is a sub-symmetric spectrum of $\operatorname{TH}(A)'[k]$, and there
is a canonical weak equivalence from the mapping cone of the inclusion
onto the symmetric spectrum $\operatorname{TH}(A/I)[k]$. Indeed, this
follows from the proof of Lemma~\ref{hypercubelemma} below. The face
maps $d_r \colon \operatorname{TH}(A)'[k] \to
\operatorname{TH}(A)'[k-1]$, $0 \leqslant r \leqslant k$, restrict
to face maps $d_r \colon \operatorname{TH}(A,I)[k] \to
\operatorname{TH}(A,I)[k-1]$, $0 \leqslant r \leqslant k$, and we define
$$\operatorname{TH}(A,I) =
\big{\|} [k] \mapsto \operatorname{TH}(A,I)[k] \big{\|}$$
to be the realization. Then from the canonical weak equivalences of
$\operatorname{TH}(A)$ and $\operatorname{TH}(A)'$ and of
$\operatorname{TH}(A/I)$ and the mapping cone of the canonical
inclusion of $\operatorname{TH}(A,I)$ in $\operatorname{TH}(A)'$ we
obtain  the following distinguished triangle in the stable homotopy
category:
$$\operatorname{TH}(A,I) \to \operatorname{TH}(A) \to
\operatorname{TH}(A/I) \xrightarrow{\partial}
\Sigma \operatorname{TH}(A,I).$$
This distinguished triangle is canonically isomorphic to the
distinguished triangle underlying the distinguished triangle of
$\mathbb{T}$-spectra of~\cite[Prop.~10.2]{bokstedtmadsen}.

The category $\mathcal{P}'[k]$ admits a finite filtration
$$\emptyset \subset \operatorname{Fil}_0\mathcal{P}'[k] \subset
\operatorname{Fil}_1\mathcal{P}'[k] \subset \dots \subset
\operatorname{Fil}_{k-1}\mathcal{P}'[k] \subset
\operatorname{Fil}_k\mathcal{P}'[k] = \mathcal{P}'[k]$$
with $\operatorname{Fil}_s\mathcal{P}'[k]$ the full sub-category with
objects the tuples $(\epsilon_0,\dots,\epsilon_k)$ such that at most
$s$ of $\epsilon_0,\dots,\epsilon_k$ are equal to $0$. We obtain an
induced filtration
$$* \subset \operatorname{Fil}_0\operatorname{TH}(A,I)[k] \subset
\operatorname{Fil}_1\operatorname{TH}(A,I)[k] \subset \dots \subset
\operatorname{TH}(A,I)[k]$$
where $\operatorname{Fil}_s\operatorname{TH}(A,I)[k]$ is the
sub-symmetric spectrum defined to be the homotopy colimit of the
restriction of the functor $E[k]$ to
$\operatorname{Fil}_s\mathcal{P}'[k]$. The functor $\mathfrak{d}_r
\colon \mathcal{P}'[k] \to \mathcal{P}'[k-1]$ takes
$\operatorname{Fil}_s\mathcal{P}'[k]$ to
$\operatorname{Fil}_s\mathcal{P}'[k-1]$, for $0 \leqslant s < k$, and
therefore, the face map $d_r \colon \operatorname{TH}(A,I)[k] \to
\operatorname{TH}(A,I)[k-1]$ takes
$\operatorname{Fil}_s\operatorname{TH}(A,I)[k]$ to
$\operatorname{Fil}_s\operatorname{TH}(A,I)[k-1]$, for
$0 \leqslant s < k$. Hence, we can form the realization
$$\operatorname{Fil}_s\operatorname{TH}(A,I) = \big{\|} [k] \mapsto
\operatorname{Fil}_s\operatorname{TH}(A,I)[k] \big{\|}$$ 
and obtain the following filtration of $\operatorname{TH}(A,I)$ by
sub-symmetric spectra:
$$* \subset \operatorname{Fil}_0\operatorname{TH}(A,I) \subset
\operatorname{Fil}_1\operatorname{TH}(A,I) \subset
\operatorname{Fil}_2\operatorname{TH}(A,I) \subset \dots \subset
\operatorname{TH}(A,I).$$
We define the filtration quotient
$\operatorname{gr}_s\operatorname{TH}(A,I)[k]$ to be the mapping cone
of the canonical inclusion of
$\operatorname{Fil}_{s-1}\operatorname{TH}(A,I)[k]$ in
$\operatorname{Fil}_s\operatorname{TH}(A,I)[k]$ and let
$$\operatorname{gr}_s\operatorname{TH}(A,I) = \big{\|} [k] \mapsto
\operatorname{gr}_s\operatorname{TH}(A,I)[k] \big{\|}.$$
Since forming mapping cone and realization commute, up to canonical
isomorphism, $\operatorname{gr}_s\operatorname{TH}(A,I)$ is
canonically isomorphic to the mapping cone of the inclusion of
$\operatorname{Fil}_{s-1}\operatorname{TH}(A,I)$ in
$\operatorname{Fil}_s\operatorname{TH}(A,I)$. We evaluate the homotopy
type of $\operatorname{gr}_s\operatorname{TH}(A,I)$ and begin by
evaluating that of $\operatorname{gr}_s\operatorname{TH}(A,I)[k]$.

Let $E^1 \to E^0$ and $E^0 \to E^{-1}$ be the maps of unital symmetric
$A$-$A$-bimodule spectra associated with the canonical inclusion of
$I$ in $A$ and the canonical projection of $A$ onto $A/I$.
We call these maps the canonical maps. We say that a tuple $\epsilon =
(\epsilon_0,\dots,\epsilon_k)$ of integers $\epsilon_i \in \{1, 0,
-1\}$ is a \emph{weight} of rank $k$ and write $|\epsilon|_v$ for the
number of $0 \leqslant i \leqslant k$ with $\epsilon_i = v$. We have
the following analog of the hyper-cube
lemma~\cite[Lemma~A.1]{wodzicki}:

\begin{lem}\label{hypercubelemma}For every $0 \leqslant s \leqslant k$, there
is a canonical weak equivalence
$$\bar{\Phi}_s \colon \operatorname{gr}_s\operatorname{TH}(A,I)[k]
\xrightarrow{\sim}
\bigvee E^{\epsilon_0} \wedge \dots \wedge E^{\epsilon_k}$$
where, on the right-hand side, the wedge sum ranges over all weights
$\epsilon$ of rank $k$ with $|\epsilon|_0 = 0$ and $|\epsilon|_{-1} =
s$.
\end{lem}

\begin{proof}There is a partial ordering on the set of weights
of rank $k$ where $\epsilon \preccurlyeq \epsilon'$ if $\epsilon_i
\geqslant \epsilon_i'$, for all $0 \leqslant i \leqslant k$. We first
define a map of symmetric spectra
$$\Phi_s \colon \operatorname{Fil}_s\operatorname{TH}(A,I)[k] =
\operatornamewithlimits{hocolim}_{\operatorname{Fil}_s\mathcal{P}'[k]}
E^{\epsilon_0} \wedge \dots \wedge E^{\epsilon_k} \to
\prod E^{\epsilon_0} \wedge \dots \wedge E^{\epsilon_k}$$
where the product on the right-hand side is indexed by the set of
weights of rank $k$ with $|\epsilon|_0 = 0$ and $|\epsilon|_{-1} =
s$. We recall from~\cite[XII.2.2]{bousfieldkan} that to give the map
$\Phi_s$ is equivalent to giving, for each weight $\epsilon$ of rank
$k$ with $|\epsilon|_0 = 0$ and $|\epsilon|_{-1} = s$, a compatible
family of maps
$$\Phi_{s,\epsilon,\epsilon'} \colon (E^{\epsilon_0'} \wedge \dots
\wedge E^{\epsilon_k'}) \wedge |
\epsilon'/\operatorname{Fil}_s\mathcal{P}'[k] |_+ \to
E^{\epsilon_0} \wedge \dots \wedge E^{\epsilon_k},$$ 
where $\epsilon'$ ranges over all objects of
$\operatorname{Fil}_s\mathcal{P}'[k]$, and where the second smash
factor on the left-hand side is the geometric realization of the nerve
of the under-category $\epsilon' /
\operatorname{Fil}_s\mathcal{P}'[k]$. Combatibility means that if
$\epsilon' \preccurlyeq \epsilon''$, then the following diagram
commutes:
$$\xymatrix{
{ (E^{\epsilon_0'} \wedge \dots \wedge E^{\epsilon_k'}) \wedge
| \epsilon'/\operatorname{Fil}_s\mathcal{P}'[k] |_+ }
\ar[r]^(.65){\Phi_{s,\epsilon,\epsilon'}}  &
{ E^{\epsilon_0} \wedge \dots \wedge E^{\epsilon_k} } \ar@{=}[dd] \cr
{ (E^{\epsilon_0'} \wedge \dots \wedge E^{\epsilon_k'}) \wedge
| \epsilon''/\operatorname{Fil}_s\mathcal{P}'[k] |_+ }
\ar[d]^{ \operatorname{can} \wedge
| \epsilon''/\operatorname{Fil}_s\mathcal{P}'[k] |_+ }
\ar[u]_{(E^{\epsilon_0'} \wedge \dots \wedge E^{\epsilon_k'}) \wedge
  \operatorname{can} } &
{ } \cr
{ (E^{\epsilon_0''} \wedge \dots \wedge E^{\epsilon_k''}) \wedge
| \epsilon''/\operatorname{Fil}_s\mathcal{P}'[k] |_+ }
\ar[r]^(.65){\Phi_{s,\epsilon,\epsilon''}} &
{ E^{\epsilon_0} \wedge \dots \wedge E^{\epsilon_k}. } \cr
}$$
Now, there is a unique weight $\omega$ of rank $k$ with $\omega
\preccurlyeq \epsilon$ and $|\omega|_0 = s$. The category
$\omega/\operatorname{Fil}_s\mathcal{P}'[k]$ has one object and one
morphism, and we define
$$\Phi_{s,\epsilon,\omega} \colon
E^{\omega_0} \wedge \dots \wedge E^{\omega_k} \to
E^{\epsilon_0} \wedge \dots \wedge E^{\epsilon_k}$$
to be the canonical map. If $\epsilon' \neq \omega$, we define
$\Phi_{s,\epsilon,\epsilon'}$ to be the constant map. To see that
these maps are compatible, suppose that $\epsilon' \preccurlyeq
\omega$. Then either $\epsilon' = \omega$ or $|\epsilon'|_0 < s$. The
former case is trivial, and in the latter case, the composition of the
following canonical maps
$$E^{\epsilon_0'} \wedge \dots \wedge E^{\epsilon_k'} \to
E^{\omega_0} \wedge \dots \wedge E^{\omega_k} \to
E^{\epsilon_0} \wedge \dots \wedge E^{\epsilon_k}$$
is equal to the constant map. Indeed, on at least one smash factor, the
map is induced by the composition $I \to A \to A/I$, which is
zero. Hence, we obtain the desired map $\Phi_s$. It is clear from the
definition that the image of $\Phi_s$ is contained in the
corresponding wedge sum, and that the restriction of $\Phi_s$ to
$\operatorname{Fil}_{s-1}\operatorname{TH}(A,I)[k]$ is equal to the
constant map. Therefore, we obtain the map $\bar{\Phi}_s$ of the
statement.

The proof that $\bar{\Phi}_s$ is a weak equivalence is by induction on
$k \geqslant 0$. The case $k = 0$ is trivial, and the induction step
uses the following observation. Suppose that $\epsilon' \preccurlyeq
\epsilon \preccurlyeq \epsilon''$ are such that for some $0 \leqslant
v \leqslant k$, $\epsilon'(v) = 1$, $\epsilon(v) = 0$, and
$\epsilon''(v) = -1$, and such that for all $0 \leqslant i \leqslant
k$ with $i \neq v$, $\epsilon'(i) = \epsilon(i) = \epsilon''(i)$. Then
the composition of the canonical maps
$$E^{\epsilon_0'} \wedge \dots \wedge E^{\epsilon_k'} \to
E^{\epsilon_0} \wedge \dots \wedge E^{\epsilon_k} \to
E^{\epsilon_0''} \wedge \dots \wedge E^{\epsilon_k''}$$
is the constant map, and the induced map from the mapping cone of the
left-hand map to the right-hand term is a weak equivalence. 
\end{proof}

We define maps of symmetric spectra
$$d_{r,s} \colon
\bigvee E^{\epsilon_0} \wedge \dots \wedge E^{\epsilon_k} \to
\bigvee E^{\epsilon_0}\wedge \dots \wedge E^{\epsilon_{k-1}},
\hskip5mm 0 \leqslant r \leqslant k,$$
where the wedge sum on the left and right-hand side range over all
weights $\epsilon$ of rank $k$ and $k-1$, respectively, with
$|\epsilon|_0 = 0$ and $|\epsilon|_{-1} = s$. Let $\epsilon$ be a
weight of rank $k$ with $|\epsilon|_0 = 0$ and $|\epsilon|_{-1} = s$,
and let $0 \leqslant r \leqslant k$. We define a weight
$\partial_r(\epsilon)$ of rank $k-1$ with $|\partial_r(\epsilon)|_0 =
0$ and with $|\partial_r(\epsilon)|_{-1} \leqslant s$ by
$$\partial_r(\epsilon)_i = \begin{cases}
\epsilon_i & \text{for $0 \leqslant i < r$,} \cr
\max\{\epsilon_r,\epsilon_{r+1}\} & \text{for $i = r$,} \cr
\epsilon_{i+1} & \text{for $r < i \leqslant k-1$,} \cr
\end{cases}$$
if $0 \leqslant r < k$, and by
$$\partial_k(\epsilon)_i = \begin{cases}
\max\{\epsilon_k,\epsilon_0\} & \text{for $i = 0$}, \cr
\epsilon_i & \text{for $0 < i \leqslant k-1$,} \cr
\end{cases}$$
if $r = k$. Then we define the component $d_{r,s,\epsilon}$ of the map
$d_{r,s}$ to be the constant map, if $|\partial_r(\epsilon)|_{-1} < s$,  
and to be the composition
$$E^{\epsilon_0} \wedge \dots \wedge E^{\epsilon_k}
\xrightarrow{\delta_r} E^{\partial_r(\epsilon)_0} \wedge \dots \wedge
E^{\partial_r(\epsilon)_{k-1}} \to
\bigvee_{\epsilon' \in \mathcal{P}(k-1,s)} E^{\epsilon_0'} \wedge
\dots \wedge E^{\epsilon_{k-1}'},$$
if $|\partial_r(\epsilon)|_{-1}  = s$. Here the right-hand map is the
canonical inclusion. The next result is an immediate consequence
of the construction of the map $\bar{\Phi}_s$:

\begin{addendum}\label{hypercubeaddendum}For all $0 \leqslant s
\leqslant k$ and $0 \leqslant r \leqslant k$, the diagram
$$\xymatrix{
{\operatorname{gr}_s\operatorname{TH}(A,I)[k]}
\ar[r]^(.46){\bar{\Phi}_s} \ar[d]^{d_r} & 
{\displaystyle{\bigvee E^{\epsilon_0} \wedge \dots \wedge
E^{\epsilon_k}}} \ar[d]^{d_{r,s}} \cr 
{\operatorname{gr}_s\operatorname{TH}(A,I)[k-1]}
\ar[r]^(.48){\bar{\Phi}_s} & 
{\displaystyle{\bigvee E^{\epsilon_0} \wedge \dots \wedge
E^{\epsilon_{k-1}}}} \cr
}$$
commutes. Here the wedge sums on the right-hand side range over all
weights $\epsilon$ of rank $k$ and $k-1$, respectively, with
$|\epsilon|_0 = 0$ and $|\epsilon|_{-1} = s$. \hfill\space\qed
\end{addendum}

We now evaluate the homotopy type of the realization
$\operatorname{gr}_s\operatorname{TH}(A,I)$. The spectra on the
right-hand side in Lemma~\ref{hypercubelemma} form a semi-simplicial
symmetric spectrum, that we denote
$\operatorname{gr}_s\operatorname{TH}(A,I)^{\dagger}[-]$, with face
maps
$$d_{r,s} \colon
\operatorname{gr}_s\operatorname{TH}(A,I)^{\dagger}[k] \to 
\operatorname{gr}_s\operatorname{TH}(A,I)^{\dagger}[k-1], \hskip 5 mm
0 \leqslant r \leqslant k,$$ 
and Addendum~\ref{hypercubeaddendum} shows that the weak equivalence
of Lemma~\ref{hypercubelemma} forms a map of semi-simplicial symmetric
spectra. The spectral sequence~(\ref{skeletonspectralsequence}) shows
that the induced map of realizations is a weak equivalence
$$\bar{\Phi}_s \colon \operatorname{gr}_s\operatorname{TH}(A,I)
\xrightarrow{\sim}
\operatorname{gr}_s\operatorname{TH}(A,I)^{\dagger}.$$
It is clear that for $s = 0$, we have
$$\operatorname{gr}_0\operatorname{TH}(A,I)^{\dagger} =
\operatorname{TH}(I).$$
We proceed to show that the higher filtration quotients are
$B(I)$-cellular, where $B(I)$ is the bar-construction of the
non-unital associative symmetric ring spectrum $I$. This is defined to
be the realization
$$B(I) = \big{\|} [k] \mapsto B(I)[k] \big{\|}$$
of the semi-simplicial symmetric spectrum with $B(I)[k] = I^{\wedge
  k}$, if $k > 0$, with $B(I)[0] = *$, and with face maps defined by
$$d_r = \begin{cases}
I^{\wedge (r-1)} \wedge \mu \wedge I^{\wedge (k-r-1)} & 0 < r < k \cr
* & \text{$r = 0$ or $r = k$.} \cr
\end{cases}$$
To this end, for $s > 0$, we further decompose
$$\operatorname{gr}_s\operatorname{TH}(A,I)^{\dagger} =
\operatorname{gr}_s\operatorname{TH}(A,I)_0^{\dagger} \vee
\operatorname{gr}_s\operatorname{TH}(A,I)_1^{\dagger}$$
where the first term on the right-hand side is the realization of the
sub-semi-simplicial symmetric spectrum that consists of the single
wedge summand with weight $\epsilon = (-1,\dots,-1,1,\dots,1)$ and
the second is the realization of the sub-semi-simplicial symmetric
spectrum that consists of the remaining wedge summands.
Following~\cite[Lemma~3.3]{wodzicki} we introduce a further filtration
of the semi-simplicial symmetric spectrum
$\operatorname{gr}_s\operatorname{TH}(A,I)_1^{\dagger}[-]$. We define
$$\operatorname{Fil}_u'\operatorname{gr}_s
\operatorname{TH}(A,I)_1^{\dagger}[k] \subset 
\operatorname{gr}_s\operatorname{TH}(A,I)_1^{\dagger}[k]$$
to consist of the wedge summands indexed by the weights $\epsilon$ for
which there exists $k-u \leqslant j \leqslant k$ with $\epsilon_j =
-1$, and define $\operatorname{gr}_u'\operatorname{gr}_s
\operatorname{TH}(A,I)_1^{\dagger}[k]$ to be the filtration
quotient. After realization, we get an induced filtration of the
symmetric spectrum
$\operatorname{gr}_s\operatorname{TH}(A,I)_1^{\dagger}$ with
filtration quotients
$$\operatorname{gr}_u'\operatorname{gr}_s
\operatorname{TH}(A,I)_1^{\dagger} = \big{\|} [k] \mapsto
\operatorname{gr}_u'\operatorname{gr}_s
\operatorname{TH}(A,I)_1^{\dagger}[k] \big{\|}.$$ 
We prove that the higher filtration quotients
$\operatorname{gr}_s\operatorname{TH}(A,I)$ are $B(I)$-cellular in the
following precise sense:

\begin{lem}\label{gradedpieces}For all positive integers $s$, and all
non-negative integers $u$, there are canonical isomorphisms in the
stable homotopy category 
$$\begin{aligned}
\Sigma^{s-1} (A/I)^{\wedge s} \wedge B(I) 
{} & \xrightarrow{\sim} \operatorname{gr}_s
\operatorname{TH}(A,I)_0^{\dagger} \cr
\bigvee \Sigma^{s+u-1} B(I)^{\wedge \ell} \wedge (A/I)^{\wedge s}
\wedge I^{\wedge u} 
{} & \xrightarrow{\sim} \operatorname{gr}_u' \operatorname{gr}_s
\operatorname{TH}(A,I)_1^{\dagger} \cr 
\end{aligned}$$
where the wedge sum ranges over all sequences $(n_0,\dots,n_{\ell})$
of non-negative integers of length $\ell \geqslant 1$ such that
$n_0 + \dots + n_{\ell} = s$, and such that $n_i > 0$, for $i > 0$.
The isomorphisms are natural with respect to the pair $(A,I)$.
\end{lem} 

\begin{proof}To establish the first isomorphism, we note that the
symmetric spectrum in simplicial degree $k$,
$\operatorname{gr}_s\operatorname{TH}(A,I)_0^{\dagger}[k]$, is
trivial, for $0 \leqslant k < s$, and canonically isomorphic to
$(A/I)^{\wedge s} \wedge B(I)[k-s+1]$, for $k \geqslant s$. Under this
isomorphism, the face maps become $d_r = (A/I)^{\wedge s} \wedge d_i$,
if $r = s - 1 + i$ with $0 < i < k-s+1$, and $d_r = *$, otherwise. The
desired isomorphism
$$\alpha \colon (A/I)^{\wedge s} \wedge B(I) \wedge
\Delta^{s-1}/\partial\Delta^{s-1} \xrightarrow{\sim}
\operatorname{gr}_s\operatorname{TH}(A,I)_0^{\dagger}$$
is then induced by the maps of symmetric spectra
$$\alpha_k \colon (A/I)^{\wedge s} \wedge B(I)[k-s+1] \wedge
\Delta_+^{k-s+1} \wedge \Delta_+^{s-1} \to
\operatorname{gr}_s\operatorname{TH}(A,I)_0^{\dagger}[k] \wedge
\Delta_+^k$$
that, in turn, are induced by the maps of pointed spaces
$$\alpha_k' \colon \Delta_+^{k-s+1} \wedge \Delta_+^{s-1} \to
\Delta_+^k$$
defined by the formula
$$\alpha_k'((x_0,\dots,x_{k-s+1}),(y_0,\dots,y_s)) =
(y_1,\dots,y_s,y_0x_0,\dots,y_0x_{k-s+1}).$$
One checks that the maps $\alpha_k$, $k \geqslant 0$, factor through
the equivalence relation that defines the realization and induce the
isomorphism $\alpha$ as stated.

The proof of the second isomorphism is similar but notationally more
involved. The symmetric spectrum $\operatorname{gr}_u'
\operatorname{gr}_s \operatorname{TH}(A,I)_1^{\dagger}[k]$ for
$0 \leqslant k < s$, is trivial, and for $k \geqslant s$, is
canonically isomorphic to the wedge sum
$$\bigvee (A/I)^{\wedge n_0} \wedge B(I)[k_1] \wedge
(A/I)^{\wedge n_1}  \wedge \dots \wedge B(I)[k_{\ell}] \wedge
(A/I)^{\wedge n_{\ell}} \wedge I^{\wedge u}$$
that ranges over all pairs of sequences of non-negative integers 
$(n_0,\dots,n_{\ell})$ and $(k_1,\dots,k_{\ell})$ such that $\ell
\geqslant 1$, such that $n_i > 0$, for $i > 0$, and such that $n_0 +
\dots + n_{\ell} = s$ and $k_1 + \dots + k_{\ell} = k - u - s +1$. The
face map $d_{r,s}$ takes the summand indexed by $(n_0,\dots,n_{\ell})$
and $(k_1,\dots,k_t,\dots,k_{\ell})$ to the summand indexed by
$(n_0,\dots,n_{\ell})$ and $(k_1,\dots,k_t-1,\dots,k_{\ell})$ by the
map induced from
$$d_i \colon B(I)[k_t] \to B(I)[k_t-1],$$
if there exists $1 \leqslant t \leqslant \ell$ and $0 < i < k_t$ such
that
$$r = n_0 + k_1 + n_1 + \dots + k_{t-1} + n_{t-1} - 1 + i,$$
and by the trivial map, otherwise. In particular, the face maps
preserve the index $(n_0,\dots,n_{\ell})$. It follows that, as a 
semi-simplicial symmetric spectrum,
$$\operatorname{gr}_u'\operatorname{gr}_s
\operatorname{TH}(A,I)_1^{\dagger}[-] = \bigvee
E(n_0,\dots,n_{\ell})[-]$$
where $E(n_0,\dots,n_{\ell})[-]$ consists of all the wedge summands
above with the index $(n_0,\dots,n_{\ell})$ fixed. To finish the proof
of the lemma, we observe that there is a natural isomorphism of
symmetric spectra
$$\beta \colon B(I)^{\wedge \ell} \wedge (A/I)^{\wedge s} \wedge
I^{\wedge u} \wedge \Delta^{s+u-1}/\partial\Delta^{s+u-1}
\xrightarrow{\sim} \big{\|} [k] \mapsto E(n_0,\dots,n_{\ell})[k]
\big{\|}.$$
The map $\beta$ is defined in a manner similar to the map $\alpha$
above and is induced from the maps of pointed spaces
$$\beta_{(k_1,\dots,k_{\ell})}' \colon \Delta_+^{k_1} \wedge \dots
\wedge \Delta_+^{k_{\ell}} \wedge \Delta_+^{s+u-1} \to \Delta_+^k$$
that to the $(\ell+1)$-tuple of points with barycentric coordinates
$$((x_{1,0},\dots,x_{1,k_1}), \dots,
(x_{\ell,0},\dots,x_{\ell,k_{\ell}}), (y_0,\dots,y_{s+u-1}))$$
associates the point with barycentric coordinates
$$\begin{aligned}
(y_1,\dots,y_{n_0}, & y_0'x_{1,0}, \dots, y_0'x_{1,k_1},y_{n_0+1},
\dots, y_{n_0+n_1}, \dots \cr
\dots , & y_0'x_{\ell,0}, \dots, y_0'x_{\ell,k_{\ell}}, 
y_{n_0+\dots+n_{\ell-1}+1},\dots,y_{s+u-1}) \cr
\end{aligned}$$
where $y_0' = y_0/\ell$. This completes the proof of the lemma.
\end{proof}

We wish to evaluate the spectrum homology
$$H_q(B(I),\mathbb{F}_p) = \pi_q(B(I) \wedge \mathbb{F}_p).$$
The absence of a unit in $I$ and the resulting absence of degeneracy
maps in $B(I)[-]$ makes this a non-trivial problem. To wit, the
Eilenberg-Zilber theorem, which states that the total complex and the
diagonal complex associated with a bi-simplicial abelian group are
chain homotopy equivalent, is not valid in the absence of degeneracy
maps. Following Wodzicki~\cite[Thm.~9.5]{wodzicki} we prove the
following result:

\begin{lem}\label{tensorproduct}Let $I$ be a non-unital associative
ring, let $p$ be a prime, and suppose that the pro-abelian group
$$\{ \operatorname{Tor}_q^{\mathbb{Z} \ltimes I^m}
(\mathbb{Z},\mathbb{Z}/p) \}_{m \geqslant 1}$$
is zero, for all positive integers $q$. Then the pro-abelian group
$$\{ H_q( B(I^m), \mathbb{Z}/p) \}_{m \geqslant 1}$$
is zero, for all integers $q$.
\end{lem}

\begin{proof}We first recall the algebraic bar-construction. Let
$\Lambda$ be a commutative ring, and let $J$ be a non-unital
associative $\Lambda$-algebra. The algebraic bar-construction of $J$
is the semi-simplicial $\Lambda$-module $B^{\Lambda}(J)[-]$ that in
degree zero is $0$, in positive degree $k$ is the $k$-fold tensor
product of copies of $J$, and whose face maps are given by
$$d_r(a_1 \otimes \dots \otimes a_k) = \begin{cases}
a_1 \otimes \dots \otimes a_ra_{r+1} \otimes \dots \otimes a_k &
0 < r < k \cr 
0 & \text{$r = 0$ or $r = k$.} \cr
\end{cases}$$
We write $B_*^{\Lambda}(J)$ for the associative chain complex of
$\Lambda$-modules, where the differential is given by the alternating
sum of the face maps. 

We show that in the skeleton spectral sequence
$$E_{s,t}^1 = \{ H_t(B(I^m)[s],\mathbb{F}_p)\} \Rightarrow
\{H_{s+t}(B(I^m),\mathbb{F}_p)\},$$
the pro-abelian groups $E_{s,t}^2$ are zero, for all integers $s$
and $t$. This implies that also the pro-abelian groups
$E_{s,t}^{\infty}$ are zero, for all integers $s$ and $t$. Since
the induced filtration of the pro-abelian groups in the abutment is
finite, the statement follows. The K\"{u}nneth formula shows that
there is a canonical isomorphism of complexes of graded
pro-abelian groups
$$\{ B_*^{\mathbb{F}_p}(H_*(I^m,\mathbb{F}_p)) \} \xrightarrow{\sim}
\{ H_*(B(I^m)[-],\mathbb{F}_p) \}  = E_{*,*}^1.$$ 
Suppose first that the additive group of the ring $I$ is torsion
free. Then the canonical map of graded $\mathbb{F}_p$-algebras
$$I^m \otimes H_*(\mathbb{Z},\mathbb{F}_p) \to H_*(I^m,\mathbb{F}_p)$$
is an isomorphism. Hence, we must show that the graded pro-abelian
group
$$\{ H_s( B_*^{\mathbb{F}_p}(I^m \otimes
H_*(\mathbb{Z},\mathbb{F}_p))) \}$$
is zero, for all integers $s$. Here we appeal to (the proof
of)~\cite[Thm.~9.5]{wodzicki}: Since $H_*(\mathbb{Z},\mathbb{F}_p)$ is
unital, it suffices to prove that the pro-abelian group
$$\{ H_s(B_*^{\mathbb{F}_p}(I^m \otimes \mathbb{F}_p)) \}$$
is zero, for all integers $s$. For $s \leqslant 0$, this pro-abelian
group is zero, by definition. For $s > 0$, the canonical maps of
pro-abelian groups
$$\{ \operatorname{Tor}_s^{\mathbb{Z} \ltimes I^m}
(\mathbb{Z},\mathbb{F}_p) \} \to
\{ H_s(B_*^{\mathbb{Z}}(I^m) \otimes \mathbb{F}_p) \} \to
\{ H_s(B_*^{\mathbb{F}_p}(I^m \otimes \mathbb{F}_p)) \}$$
are isomorphisms, since $I$ is torsion
free~\cite[Lemma~1.1]{suslin4}, and the pro-abelian group on the
left-hand side is zero by assumption. This completes the proof for $I$
torsion free. The general case is proved in a similar manner by first
functorially replacing the non-unital associative $I^m$ by a
simplicial non-unital associative ring $F(I^m)[-]$ such that for all
$s \geqslant 0$, the non-unital associative ring $F(I^m)[s]$ is
torsion free as an abelian group.
\end{proof}

\begin{proof}[Proof of Thm.~\ref{mainthh}]We prove that the canonical map
$$\{ \operatorname{TH}_q(I^m,\mathbb{Z}/p) \}_{m \geqslant 1} \to
\{ \operatorname{TH}_q(A,I^m,\mathbb{Z}/p) \}_{m \geqslant 1}$$
is an isomorphism of pro-abelian groups, for all integers $q$. It
suffices by~\cite[Cor.~5.8]{panin} to show that the canonical map
$$\{ H_q(\operatorname{TH}(I^m),\mathbb{Z}/p) \}_{m \geqslant 1} \to
\{ H_q(\operatorname{TH}(A,I^m),\mathbb{Z}/p) \}_{m \geqslant 1}$$
is an isomorphism of pro-abelian groups, for all integers $q$. This
map is the edge-homomorphism in the spectral sequence of pro-abelian
groups
$$E_{s,t}^1 = \{ H_t(\operatorname{gr}_s
\operatorname{TH}(A,I^m),\mathbb{Z}/p) \} \Rightarrow \{
H_{s+t}(\operatorname{TH}(A,I^m),\mathbb{Z}/p) \}.$$
Hence, it suffices to show that for all $s > 0$, the pro-abelian
groups $E_{s,t}^1$ are zero. We have a direct sum decomposition
$$E_{s,t}^1 = H_t(\operatorname{gr}_s
\operatorname{TH}(A,I^m)_0^{\dagger},\mathbb{Z}/p)
\oplus H_t(\operatorname{gr}_s
\operatorname{TH}(A,I^m)_1^{\dagger},\mathbb{Z}/p)$$
and an additional first-quadrant spectral sequence
$${}'E_{u,v}^1 = \{ H_v(\operatorname{gr}_u'\operatorname{gr}_s
\operatorname{TH}(A,I^m)_1^{\dagger},\mathbb{Z}/p) \}
\Rightarrow \{ H_{u+v}(\operatorname{gr}_s
\operatorname{TH}(A,I^m)_1^{\dagger},\mathbb{Z}/p).$$
It follows that to prove the theorem, it will suffice to show that the
pro-abelian groups $\{H_q(\operatorname{gr}_s
\operatorname{TH}(A,I^m),\mathbb{Z}/p)_0^{\dagger}\}$ and
$\{H_q(\operatorname{gr}_u'\operatorname{gr}_s
\operatorname{TH}(A,I^m),\mathbb{Z}/p)_1^{\dagger}\}$ are zero, for
all positive integers $s$ and all non-negative integers $u$ and
$q$. But this is an immediate consequence of Lemmas~\ref{gradedpieces}
and~\ref{tensorproduct}.
\end{proof}

\section{Proof of Theorems~\ref{main}, \ref{curves},
and~\ref{grouprings}}\label{proofoftheorem}

In this section we complete the proof of Thm.~\ref{main} and prove
Thms.~\ref{curves} and~\ref{grouprings}. The following result is
similar to the proof of~\cite[Thm. 5.3]{cuntzquillen}:

\begin{lem}\label{freeisenough}Suppose that Thm.~\ref{main} holds for
every triple $(A,B,I)$, where $I$ can be embedded as a two-sided ideal
of a free unital associative ring. Then Thm.~\ref{main} holds for all
triples $(A,B,I)$.
\end{lem}

\begin{proof}Let $X$ be a functor from unital associative rings to
symmetric spectra that satisfies that for every pair of rings $(A,B)$
the canonical map
$$X(A \times B) \to X(A) \times X(B)$$
is a weak equivalence. We extend $X$ to a functor $X'$ defined on
non-unital associative rings as follows. If $I$ is a non-unital
associative ring, we define $X'(I)$ to be the homotopy fiber of the
map $X(\mathbb{Z} \ltimes I) \to X(\mathbb{Z})$ induced by the canonical
projection. Suppose that $A$ is a unital associative ring with unit
$e$. Then the map $\mathbb{Z} \ltimes A \to \mathbb{Z} \times A$ that
takes $(x,a)$ to $(x,x \cdot e + a)$ is a ring isomorphism, and hence,
we obtain a weak equivalence
$$X(\mathbb{Z} \ltimes A) \xrightarrow{\sim} X(\mathbb{Z} \times A)
\xrightarrow{\sim} X(\mathbb{Z}) \times X(A)$$
of symmetric spectra over $X(\mathbb{Z})$. These maps induce a weak
equivalence of the homotopy fiber $X'(A)$ of the map $X(\mathbb{Z}
\ltimes A) \to X(\mathbb{Z})$ and the homotopy fiber $X''(A)$ of the
canonical projection from $X(\mathbb{Z}) \times X(A)$ to
$X(\mathbb{Z})$. This shows that we have natural weak equivalences
$$X'(A) \xrightarrow{\sim} X''(A) \xleftarrow{\sim} X(A).$$
Hence, the functor $X'$ extends, up to canonical weak
equivalence, the functor $X$ to non-unital associative rings. In the
following, we will write $X$ instead of $X'$ for the extended functor.

Let $X^n(A)$ be the homotopy fiber of the cyclotomic trace map
$$K(A) \to \operatorname{TC}^n(A;p).$$
Then Thm.~\ref{main} is equivalent to the statement that for every
unital associative ring $A$ and two-sided ideal $I \subset A$, for all
integers $q$, all primes $p$, and all positive integers $v$, the
canonical map
$$\{X_q^n(I,\mathbb{Z}/p^v)\}_{n \geqslant 1} \to
\{X_q^n(A,I,\mathbb{Z}/p^v)\}_{n \geqslant 1}$$
is an isomorphism of pro-abelian groups. We let $F(A)$ denote the
free unital associative ring generated by the underlying set of
$A$. There is a canonical surjective ring homomorphism $F(A) \to
A$. We consider the following diagram of non-unital associative
rings, where the rows and columns of the underlying diagram of abelian
groups all are short-exact sequences:
$$\xymatrix{
{J(A,I)} \ar[r] \ar[d] &
{J(A)} \ar[r] \ar[d] &
{J(A/I)} \ar[d] \cr
{F(A,I)} \ar[r] \ar[d] &
{F(A)} \ar[r] \ar[d] &
{F(A/I)} \ar[d] \cr
{I} \ar[r] &
{A} \ar[r] &
{A/I.} \cr
}$$
Since the non-unital rings in the top row can be embedded as ideals of
free unital associative rings, the diagram, by hypothesis, gives rise
to the following map of long-exact sequences of pro-abelian groups:
$$\xymatrix{
{:} \ar[d] &
{:} \ar[d] \cr
{\{X_q^n(J(A,I),\mathbb{Z}/p^v)\}} \ar[r]^(.45){f_*'} \ar[d] &
{\{X_q^n(J(A),J(A,I),\mathbb{Z}/p^v)\}} \ar[d] \cr
{\{X_q^n(F(A,I),\mathbb{Z}/p^v)\}} \ar[r]^(.45){f_*} \ar[d] &
{\{X_q^n(F(A),F(A,I),\mathbb{Z}/p^v)\}} \ar[d] \cr
{\{X_q^n(I,\mathbb{Z}/p^v)\}} \ar[r]^(.45){f_*''} \ar[d] &
{\{X_q^n(A,I,\mathbb{Z}/p^v)\}} \ar[d] \cr
{:} &
{:} \cr
}$$
Finally, since the non-unital associative rings $J(A,I)$ and $F(A,I)$
can be embedded in a free unital associative ring, the maps labelled
$f_*'$ and $f_*$ are isomorphisms of pro-abelian groups. Hence, so is
the map labelled $f_*''$. This completes the proof.
\end{proof}

We recall the following result of Cuntz and
Quillen~\cite[Prop.~4.1]{cuntzquillen}. Let $k$ be a field, and let
$I \subset T$ be a two-sided ideal of a free unital associative
$k$-algebra. Then there exists a $k$-linear section 
\begin{equation}\label{alpha}
\alpha \colon I \to T \otimes_k I
\end{equation}
of the multiplication map $\mu \colon T \otimes_k I \to I$ such that the
diagram
$$\xymatrix{ {T \otimes_k I} \ar[d]^{\mu}
\ar[rr]^(.43){\operatorname{id} \otimes_k \alpha} &&
{T \otimes_k T \otimes_k I}
\ar[d]^{\mu \otimes_k \operatorname{id}} \cr 
{I} \ar[rr]^(.43){\alpha} && {T \otimes_k I}  \cr }$$
commutes. This is stated in~\emph{loc.cit.}~for $I$ a two-sided ideal
of free unital associative $\mathbb{C}$-algebra, but the proof works
for every ground field. The diagram shows, in particular, that the map
$\alpha$ takes $I^2$ to $I \otimes_k I$. We use this to prove the
following result:

\begin{prop}\label{hunital}Let $I \subset T$ be a two-sided ideal of
a free unital associative ring. Then for all primes $p$, and all
positive integer $v$ and $q$, the following pro-abelian group
is zero:
$$\{ \operatorname{Tor}_q^{\mathbb{Z} \ltimes I^m}
(\mathbb{Z},\mathbb{Z}/p^v) \}_{m \geqslant 1}.$$
\end{prop}

\begin{proof}Since $\mathbb{Z} \ltimes I^m$ is flat over $\mathbb{Z}$,
the group $\operatorname{Tor}_q^{\mathbb{Z} \ltimes I^m}
(\mathbb{Z},\mathbb{Z}/p^v)$ is canonically isomorphic to the $q$th
homology group of the bar-complex
$\smash{B_*^{\mathbb{Z}}(I^m) \otimes \mathbb{Z}/p^v}$; see the proof
of Lemma~\ref{tensorproduct}. Hence, it suffices to show that for
every $m \geqslant 1$, there exists $m' \geqslant m$ such that the
canonical map
$$B_*^{\mathbb{Z}}(I^{m'}) \otimes \mathbb{Z}/p^v \to
B_*^{\mathbb{Z}}(I^m) \otimes \mathbb{Z}/p^v$$
induces the zero map on homology in positive degrees. Since also
$I^m \subset T$ is a two-sided ideal of a free unital associative
ring, it suffices to treat the case $m = 1$. By simple induction, it
also suffices to consider the case $v = 1$.

The tensor-product $I \otimes \mathbb{F}_p$ is a non-unital
associative $\mathbb{F}_p$-algebra, and there is a canonical
isomorphism of complexes of $\mathbb{F}_p$-vector spaces
$$B_*^{\mathbb{Z}}(I) \otimes \mathbb{F}_p \xrightarrow{\sim}
B_*^{\mathbb{Z}}(I \otimes \mathbb{F}_p).$$
We let $(I \otimes \mathbb{F}_p)'$ and
$\overline{I \otimes \mathbb{F}_p}$ denote the kernel and image,
respectively, of the map $I \otimes \mathbb{F}_p \to T \otimes
\mathbb{F}_p$ induced by the inclusion. Then we have the following
extension of non-unital associative $\mathbb{F}_p$-algebras:
$$(I \otimes \mathbb{F}_p)' \to
I \otimes \mathbb{F}_p \to
\overline{I \otimes \mathbb{F}_p}.$$
The non-unital associative $\mathbb{F}_p$-algebra
$\overline{I \otimes \mathbb{F}_p}$ is a two-sided ideal of the free
unital associative $\mathbb{F}_p$-algebra $T \otimes \mathbb{F}_p$. Let 
$$\alpha \colon \overline{I \otimes \mathbb{F}_p} \to
(T \otimes \mathbb{F}_p) \otimes (\overline{I \otimes \mathbb{F}_p})$$
be the section of the multiplication from~(\ref{alpha}). Then the
formula
$$s(x_1 \otimes \dots \otimes x_n) = (-1)^n x_1 \otimes \dots \otimes
x_{n-1} \otimes \alpha(x_n)$$
defines a null-homotopy of the chain map
$$B_*^{\mathbb{Z}}(\overline{I^2 \otimes \mathbb{F}_p}) =
B_*^{\mathbb{Z}}((\overline{I \otimes \mathbb{F}_p})^2) \to
B_*^{\mathbb{Z}}(\overline{I \otimes \mathbb{F}_p})$$
induced by the canonical inclusion.

We next let $C_*(I \otimes \mathbb{F}_p)$ be the kernel complex of the
chain-map
$$B_*^{\mathbb{Z}}(I \otimes \mathbb{F}_p) \to
B_*^{\mathbb{Z}}(\overline{I \otimes \mathbb{F}_p})$$
induced by the canonical projection. The $\mathbb{F}_p$-vector space
$C_n(I \otimes \mathbb{F}_p)$ is generated by the tensors $x_1 \otimes
\dots \otimes x_n$, where at least one factor is in $(I \otimes
\mathbb{F}_p)'$. We claim that the chain-map
$$C_*(I^2 \otimes \mathbb{F}_p) \to C_*(I \otimes \mathbb{F}_p)$$
induced by the canonical inclusion is equal to zero. This follows
immediately from the statement that the map induced by the canonical
inclusion
$$\xymatrix{
{ (I^2 \otimes \mathbb{F}_p)' } \ar[r] \ar@{=}[d] &
{ (I \otimes \mathbb{F}_p)' } \ar@{=}[d] \cr
{ (I^2 \cap pT) / pI^2 } \ar[r] &
{ (I \cap pT) / pI } \cr }$$
is zero.
Hence, it suffices to show that $I^2 \cap pT \subset pI$, or
equivalently, that if $x \in I$ and $y \in I$, and if $xy \in pT$,
then $xy \in pI$. Now, if $xy \in pT$, then $x \in pT$ or $y \in pT$,
since a free unital associative $\mathbb{F}_p$-algebra does not have
left or right zero-divisors~\cite[p.~32]{cohn}.
If $x \in pT$, we write $x = px'$ with $x' \in T$. Then
$xy = px'y$ and since $I$ is a left ideal, $x'y \in I$. Similarly, if
$y \in pT$, we write $y = py'$ with $y' \in T$. Then $xy = pxy'$, and
since $I$ is a right ideal, $xy' \in I$. In either case, $xy \in pI$
as desired.

Finally, a diagram chase shows that the chain map
$$B_*^{\mathbb{Z}}(I^4 \otimes \mathbb{F}_p) \to B_*^{\mathbb{Z}}(I
\otimes \mathbb{F}_p)$$
induces the zero map on homology in positive degrees.
\end{proof}

\begin{lem}\label{nil}Let $f \colon A \to B$ be a map of unital
associative rings and let $I \subset A$ be a two-sided ideal such
that $f \colon I \to f(I)$ is an isomorphism onto an ideal of
$B$. Let $p$ be a prime and let $v$ be a positive integer. 
Suppose that the cyclotomic trace induces an isomorphism of
pro-abelian groups
$$\{K_q(A,B,I^m,\mathbb{Z}/p^v)\}_{m \geqslant 1} \xrightarrow{\sim}
\{\operatorname{TC}_q^n(A,B,I^m;p,\mathbb{Z}/p^v)\}_{m,n \geqslant 1},$$
for all integers $q$. Then the map induced by the cyclotomic trace
$$K_q(A,B,I,\mathbb{Z}/p^v) \to
\{\operatorname{TC}_q^n(A,B,I;p,\mathbb{Z}/p^v)\}_{n \geqslant 1}$$
is an isomorphism of pro-abelian groups, for all integers $q$.
\end{lem}

\begin{proof}We note that for all positive integers $m$ and $n$, we
have a commutative diagram of spectra in which the columns are
distinguished triangles.
$$\xymatrix{
{ K(A,B,I^m) } \ar[r] \ar[d] &
{ \operatorname{TC}^n(A,B,I^m;p) } \ar[d] \cr
{ K(A,B,I) } \ar[r] \ar[d] &
{ \operatorname{TC}^n(A,B,I;p) } \ar[d] \cr
{ K(A/I^m,B/I^m,I/I^m) } \ar[r] \ar[d] &
{ \operatorname{TC}^n(A/I^m,B/I^m,I/I^m;p) } \ar[d] \cr
{ \Sigma K(A,B,I^m) } \ar[r] &
{ \Sigma \operatorname{TC}^n(A,B,I^m;p) }. \cr
}$$
The induced diagram of homotopy groups with
$\mathbb{Z}/p^v$-coefficients gives rise, as $m$ and $n$ vary, to a
map of long-exact sequences of pro-abelian groups. The map of
pro-abelian groups induced by the top horizontal map in the diagram
above is an isomorphism by assumption, and we wish to show that the
map of pro-abelian groups induced by the second horizontal map in the
diagram above is an isomorphism. By a five-lemma argument this is
equivalent to showing that the map of pro-abelian groups induced from
the third horizontal map in the diagram above is an isomorphism. Now,
it follows from~\cite{mccarthy1} and from~\cite[Thm.~2.2.1]{gh3} that
for $m \geqslant 1$ fixed, the cyclotomic trace induces an isomorphism
of pro-abelian groups
$$K_q(A/I^m,B/I^m,I/I^m,\mathbb{Z}/p^v) \xrightarrow{\sim}
\{\operatorname{TC}_q^n(A/I^m,B/I^m,I/I^m;
p,\mathbb{Z}/p^v)\}_{n \geqslant 1}.$$
But then the map of pro-abelian groups obtained by also allowing
$m \geqslant 1$ to vary is an isomorphism. This completes the proof.
\end{proof}

\begin{proof}[Proof of Thm.~\ref{main}]By Lemma~\ref{freeisenough}, we
can assume that the non-unital associative ring $I$ can be embedded as
a two-sided ideal $I \subset T$ of a free unital associative
ring. Moreover, by Lemma~\ref{nil} it suffices to show that for all
integers $q$, all primes $p$, and all positive integers $v$, the
cyclotomic trace induces an isomorphism of pro-abelian groups
$$\{K_q(A,B,I^m,\mathbb{Z}/p^v)\}_{m \geqslant 1} \xrightarrow{\sim}
\{\operatorname{TC}_q^n(A,B,I^m;p,\mathbb{Z}/p^v)\}_{m,n \geqslant 1}.$$
But Prop.~\ref{hunital} and Thms.~\ref{ktheory} and~\ref{maintc} show
that both pro-abelian groups are zero.
\end{proof}

Finally, we prove Thms.~\ref{curves} and~\ref{grouprings} of
the introduction.

\begin{proof}[Proof of Thm.~\ref{curves}]Let $X \to \operatorname{Spec} k$
be a curve over a field $k$, that is, a separated morphism of finite
type with $X$ of dimension one. We wish to show that for all
$q \geqslant r$, the cyclotomic trace
$$K_q(X,\mathbb{Z}/p^v) \to
\{\operatorname{TC}_q^n(X;p,\mathbb{Z}/p^v)\}_{n \geqslant 1}$$
is an isomorphism of pro-abelian groups. We first reduce to the case
where $X$ is affine. A covering $\{U_i\}$ of $X$ by affine open
subschemes gives rise to spectral sequences of pro-abelian groups
$$E_{s,t}^2 = \check{H}^{-s}(\{U_i\},K_t(-,\mathbb{Z}/p^v))
\Rightarrow K_{s+t}(X,\mathbb{Z}/p^v)$$
and
$$E_{s,t}^2 = \{\check{H}^{-s}(\{U_i\},
\operatorname{TC}_t^n(-;p,\mathbb{Z}/p^v))\}_{n \geqslant 1}
\Rightarrow
\{\operatorname{TC}_{s+t}^n(X;p,\mathbb{Z}/p^v)\}_{n \geqslant 1}.$$
Indeed, for $K$-theory, this follows
from~\cite[Thm.~10.3]{thomasontrobaugh}, and for topological cyclic
homology, this follows directly from the definition of the
topological cyclic homology of a
scheme~\cite[Def.~3.2.3]{gh}. Assuming the theorem for affine schemes,
the cyclotomic trace induces a map of spectral sequences which is an
isomorphism on $E_{s,t}^2$ for $t \geqslant r$. But then the
cyclotomic trace induces an isomorphism of the pro-abelian groups
in the abutments in degrees $q \geqslant r$.

So assume that $X$ is affine. By an argument similar to the proof of
Lemma~\ref{nil} above, we can further assume that $X$ is reduced. By
induction on the number of irreducible components in $X$, we can
further reduce to the case where $X$ is irreducible. Indeed, let $Z$
and $Z'$ be two proper closed subschemes of $X$ and suppose that $Z$
and $Z'$ cover $X$, or equivalently, that the corresponding ideals $I$
and $I'$ of the coordinate ring $A$ intersect trivially. Then it
follows from Thm.~\ref{main} that the cyclotomic trace induces an
isomorphism of pro-abelian groups
$$K_q(A,A/I',I,\mathbb{Z}/p^v) \xrightarrow{\sim}
\{\operatorname{TC}_q^n(A,A/I',I;p,\mathbb{Z}/p^v)\}_{n \geqslant 1},$$
for all integers $q$. This proves the induction step.

Finally, let $X \to \operatorname{Spec} k$ be a reduced and
irreducible affine curve. The coordinate ring $A$ is an integral
domain, and we let $B$ be the integral closure. We claim that the
conductor ideal $I = \{a \in A \mid aB \subset A\}$, that is the
largest ideal of $A$ which is also an ideal of $B$, is
non-zero. Indeed, the extension $A \subset B$ is finite, so we can
write $B = b_1 A + \dots +b_m A$. Since $A$ and $B$ have the same
quotient field, there are non-zero elements $a_i \in A$ such that
$a_ib_i \in A$, $1 \leqslant i \leqslant m$. But then the product
$a_1\dots a_m$ is in $I$ which therefore is non-zero as claimed. The
ring $B$ is a regular $k$-algebra, and since $A$ and $B$ are of
dimension one, the quotients $A/I$ and $B/I$ are finite
$k$-algebras. Hence, we know from \cite[Thm.~4.2.2]{gh} and
from~\cite{mccarthy1} and~\cite[Thm.~2.2.1]{gh3} that the cyclotomic
trace induces an isomorphism of pro-abelian groups
$$K_q(R,\mathbb{Z}/p^v) \xrightarrow{\sim}
\{\operatorname{TC}_q^n(R;p,\mathbb{Z}/p^v)\}_{n \geqslant 1},$$
for $q \geqslant r$, if $R$ is one of $B$, $A/I$, and
$B/I$. But then Thm.~\ref{main} shows that this map is an isomorphism,
for $q \geqslant r$, also if $R = A$.
\end{proof}

\begin{proof}[Proof of Thm.~\ref{grouprings}]The ideal $I = |G| \cdot
\mathfrak{M}$ is a common ideal of $A = \mathbb{Z}[G]$ and $B =
\mathfrak{M}$~\cite[Cor.~XI.1.2]{bass}, and hence, we can apply
Thm.~\ref{main}. It follows that the exactness of the sequence in the
statement of Thm.~\ref{grouprings} is equivalent to the exactness of
the corresponding sequence with $\mathbb{Z}[G]$ and $\mathfrak{M}$
replaced by $\mathbb{Z}[G]/I$ and $\mathfrak{M}/I$. The latter are
both artinian rings whose quotient by the Jacobson radical is a
semi-simple $\mathbb{F}_p$-algebra. Hence, the exactness of the latter
sequence follows from~\cite[Thm.~5.1(i)]{hm} and
from~\cite[Thm.~2.2.1]{gh3}.
\end{proof}

\providecommand{\bysame}{\leavevmode\hbox to3em{\hrulefill}\thinspace}
\providecommand{\MR}{\relax\ifhmode\unskip\space\fi MR }
\providecommand{\MRhref}[2]{%
  \href{http://www.ams.org/mathscinet-getitem?mr=#1}{#2}
}
\providecommand{\href}[2]{#2}

\end{document}